\documentclass[8 pt]{amsproc}
\usepackage{amssymb}
\oddsidemargin=-1cm\evensidemargin=-1cm
\begin{document}
\numberwithin{equation}{section}

\def\1#1{\overline{#1}}
\def\2#1{\widetilde{#1}}
\def\3#1{\widehat{#1}}
\def\4#1{\mathbb{#1}}
\def\5#1{\frak{#1}}
\def\6#1{{\mathcal{#1}}}

\def\C{{\4C}}
\def\R{{\4R}}
\def\n{{\4n}}
\def\Z{{\4Z}}

\author{Valentin Burcea  }
\title{Special Classes of  CR Singularities I}
 \begin{abstract}It is the sequel of the previous work about CR Singularities. There are proven few new analogues of the well-known Theorem of Moser for Special Classes of Real-Analytic CR Singular Real Submanifolds in Complex Spaces. The applied method depends on an iterative normalization procedure. \end{abstract}
\address{V. Burcea: INDEPENDENT}
\email{valentin@maths.tcd.ie}
\thanks{ THANKS ($50\%$) to   Science Foundation Ireland Grant 10/RFP/MTH 2878 for ENORMOUS FUNDING while I was working in Trinity College Dublin, Ireland}
\thanks{THANKS ($50\%$) to Capes for GENEROUS FUNDING while I was working in The Federal University of Santa Catarina, Brazil} 
\thanks{Emphasizing that the reference \cite{bu1} was fully supported by  Science Foundation Ireland Grant 06/RFP/MAT 018}
 
\dedicatory{To The Memory of Professor Steriu Ianu\c{s}}
  \maketitle 
  
 \def\Label#1{\label{#1}{\bf (#1)}~}


\def\cn{{\C^n}}
\def\cnn{{\C^{n'}}}
\def\ocn{\2{\C^n}}
\def\ocnn{\2{\C^{n'}}}


\def\dist{{\rm dist}}
\def\const{{\rm const}}
\def\rk{{\rm rank\,}}
\def\id{{\sf id}}
\def\tr{{\bf tr\,}}
\def\aut{{\sf aut}}
\def\Aut{{\sf Aut}}
\def\CR{{\rm CR}}
\def\GL{{\sf GL}}
\def\Re{{\sf Re}\,}
\def\Im{{\sf Im}\,}
\def\span{\text{\rm span}}
\def\Diff{{\sf Diff}}

\def\codim{{\rm codim}}
\def\crd{\dim_{{\rm CR}}}
\def\crc{{\rm codim_{CR}}}

\def\phi{\varphi}
\def\eps{\varepsilon}
\def\d{\partial}
\def\a{\alpha}
\def\b{\beta}
\def\g{\gamma}
\def\G{\Gamma}
\def\D{\Delta}
\def\Om{\Omega}
\def\k{\kappa}
\def\l{\lambda}
\def\L{\Lambda}
\def\z{{\bar z}}
\def\w{{\bar w}}
\def\Z{{\1Z}}
\def\t{\tau}
\def\th{\theta}

\emergencystretch15pt \frenchspacing

\newtheorem{Thm}{Theorem}[section]
\newtheorem{Cor}[Thm]{Corollary}
\newtheorem{Pro}[Thm]{Proposition}
\newtheorem{Lem}[Thm]{Lemma}

\theoremstyle{definition}\newtheorem{Def}[Thm]{Definition}

\theoremstyle{remark}
\newtheorem{Rem}[Thm]{Remark}
\newtheorem{Exa}[Thm]{Example}
\newtheorem{Exs}[Thm]{Examples}

\def\bl{\begin{Lem}}
\def\el{\end{Lem}}
\def\bp{\begin{Pro}}
\def\ep{\end{Pro}}
\def\bt{\begin{Thm}}
\def\et{\end{Thm}}
\def\bc{\begin{Cor}}
\def\ec{\end{Cor}}
\def\bd{\begin{Def}}
\def\ed{\end{Def}}
\def\br{\begin{Rem}}
\def\er{\end{Rem}}
\def\be{\begin{Exa}}
\def\ee{\end{Exa}}
\def\bpf{\begin{proof}}
\def\epf{\end{proof}}
\def\ben{\begin{enumerate}}
\def\een{\end{enumerate}}
\def\beq{\begin{equation}}
\def\eeq{\end{equation}}

 \section{Introduction and Main Results}
In this paper, there are  considered   Special Classes of Real-Formal C.-R. Singular Submanifolds  in Complex Spaces, which are formally (holomorphically) equivalent\cite{ko2},\cite{ko3}, in order to study the (local) Equivalence Problem in Complex Analysis.  This means the existence of a formal holomorphic equivalence\cite{mi1},\cite{mi2},\cite{memiza},\cite{su}  between such Real  Submanifolds.  Concerning the non-triviality of the convergence\cite{mi1},\cite{mi2} of such formal equivalence,  Moser-Webster\cite{mowe}, Gong\cite{go2} and Kossovskiy-Shafikov\cite{kosh} studied the existence   of two Real-Snalytic Submanifolds in Complex Spaces, which are only formally (holomorphically) equivalent.  Their concrete examples show the existence of a beautiful phenomenon: the  analiticity of the Real Submanifolds does not imply the analyticity of their formal equivalences\cite{go2},\cite{mowe},\cite{kosh}.
 
In this paper, it is approached the (local)  Equivalence Problem being motivated by Moser's Theorem\cite{mo} and its generalizations \cite{bu3},\cite{bu4},\cite{val8},\cite{go1}, \cite{huyi1}. It has been an importance source of inspiration for exploring the (local) Equivalence Problem. In particular, the author\cite{bu3} generalized Moser's Theorem\cite{mo}  considering a very Special Class of Real-Analytic Submanifolds in Complex Spaces: 

Let $\left(z_{11},\dots,z_{1N},\dots,z_{m1},\dots,z_{mN},w_{11},\dots,w_{mm}\right)$ be the coordinates in $\mathbb{C}^{mN+m^{2}}$, where $m,N\in\mathbb{N}^{\star}$. Introducing the following  matrices $$W=\left\{w_{ij}\right\}_{1\leq i,j\leq m},\quad Z=\left\{z_{ij}\right\}_{1\leq i\leq m,\hspace{0.1 cm}1\leq j \leq N},$$ 
there was considered the following class of Real-Analytic Submanifolds 
\begin{equation} W=Z\overline{Z}^{t}+\rm{O}\left(3\right),\label{5veve}\end{equation}
Then, the author\cite{bu3} showed that (\ref{5veve}) is  (bi)holomorphically equivalent to the following  model
\begin{equation} W=Z\overline{Z}^{t},\label{5veve1}\end{equation}
 if and only  if (\ref{5veve}) is formally equivalent to it.
 
This result\cite{bu3} originates from Moser\cite{mo},   when $m=1$ and $N=1$, and respectively from  Huang-Yin\cite{huyi1},  when $m=1$ and $N\in\mathbb{N}^{\star}$. Other  similar results have been established by Gong\cite{go1} in $\mathbb{C}^{2}$, and then later by the author\cite{bu5},\cite{val8}   in $\mathbb{C}^{N+1}$. Also, Mir\cite{mi1},\cite{mi2} discovered   other similar convergence results in the C.-R. Setting.

Clearly, any biholomorphic equivalence is a formal equivalence, but it  is not true conversely. In particular, it is not trivial to show  the convergence of a formal (holomorphic) mapping\cite{mi1},\cite{mi2}. The author\cite{bu3}  approached this convergence problem by applying   the arguments of   rapid convergence from Moser\cite{mo} and Huang-Yin\cite{huyi1}, after it was formally established a partial normal form using Generalized Fischer Decompositions\cite{sh}. Then, the author\cite{bu3} showed the convergence of a formal equivalence which was normalized by suitable compositions using  formal automorphisms of the   model (\ref{5veve1}).  

In this paper, we study   the (local) Equivalence Problem from the following points of view: 
\begin{itemize}
\item the equivalence of  formal embeddings defined between Model-Manifolds;
\item the equivalence problem of   Real-Formal Submanifolds, which are formally embedded in   Model-Manifolds;
\item  the equivalence of   formal embeddings defined  between  Real-Formal Submanifolds, which are formally embedded in  Model-Manifolds.
\end{itemize}

It is   important to  emphasize that these (formal) equivalences are considered up to suitable compositions with formal automorphisms of the corresponding Model-Manifolds, defining  consequently  classes of equivalence as in  \cite{BH},\cite{BEH},\cite{Far},\cite{H},\cite{huang1},\cite{huang2},\cite{HJ},\cite{HJX},\cite{HJY}. The   discovered gapping phenomenas  \cite{BH},\cite{BEH},\cite{H},\cite{huang1},\cite{huang2},\cite{HJX},\cite{HJY}    brought novelty in Complex Analysis   and  an important source of inspiration motivating the present research studies. 

The computations are   based on formal expansions in the local defining equations. In particular, we consider Special Classes of Real-Formal CR Singular Submanifolds in Complex Spaces\cite{bu11}, and respectively C.-R. Singular Model-Manifolds situated in Complex Spaces of possibly different dimensions. We recall that $p=0$ is called C.-R. Singularity for the Real Submanifold $M\subset\mathbb{C}^{N+1}$ if $p=0$ is a jumping point for the   mapping $M\ni q\mapsto \dim_{\mathbb{C}}T^{c}_{q}M$ defined near $p=0$.  We denoted by  $T^{c}_{q}M$    the complex tangent space of $M$ at $q$.  The reader is guided to  \cite{bu11} and \cite{gole} for studies regarding  classes of C.-R. Singularities in Codimension $2$.

In this paper, we focus on C.-R. Singular Real-Formal Submanifolds in Complex Spaces and of Codimension different from $2$,   with respect to the following situation. Let $\left(w,z\right)=\left(w_{1},w_{2},\dots,w_{N};z_{1},z_{2},\dots,z_{N}\right)$ be the coordinates in $\mathbb{C}^{2N}$. Respectively, let $\left(w',z'\right)=\left({w'}_{1},{w'}_{2},\dots,{w'}_{N};{z'}_{1},{z'}_{2},\dots,{z'}_{N}\right)$ be the coordinates in $\mathbb{C}^{2N'}$. The first result is the following:

\bt\label{ta} Let 
    model-manifolds $\mathcal{M}\subset\mathbb{C}^{2N}$ and $\mathcal{M}'\subset\mathbb{C}^{2N'}$  defined as follows
 \begin{equation}\mathcal{M}_{N}:\hspace{0.1 cm}\left\{\begin{split}&w_{1}=z_{1}\overline{z}_{1}+\lambda_{1}\left(z_{1}^{2}+\overline{z}_{1}^{2}\right) ,\\& w_{2}=z_{2}\overline{z}_{2}+\lambda_{2}\left(z_{2}^{2}+\overline{z}_{2}^{2}\right) ,\\&\quad\quad\vdots \\& w_{N}=z_{N}\overline{z}_{N}+\lambda_{N}\left(z_{N}^{2}+\overline{z}_{N}^{2}\right) ,\end{split}\right.\quad\mathcal{M}'_{N'}:\hspace{0.1 cm}\left\{\begin{split}&{w'}_{1}={z'}_{1}\overline{{z'}}_{1}+\lambda'_{1}\left({z'}_{1}^{2}+\overline{{z'}}_{1}^{2}\right) ,\\& {w'}_{2}={z'}_{2}\overline{{z'}}_{2}+\lambda'_{2}\left({z'}_{2}^{2}+\overline{{z'}}_{2}^{2}\right) ,\\&\quad\quad\vdots \\& {w'}_{N'}={z'}_{N'}\overline{{z'}}_{N'}+\lambda'_{N'}\left({z'}_{N'}^{2}+\overline{{z'}}_{N'}^{2}\right) ,\end{split}\right.\label{emo}\end{equation}
where it is assumed the following 
 \begin{equation}
\lambda_{1}, \lambda_{2}, \dots,\lambda_{N} \in \left(0,\frac{1}{2}\right),\quad \lambda'_{1}, \lambda'_{2}, \dots,\lambda'_{N'} >0. \label{lambda}
 \end{equation}
 
 Then, if $F: \mathcal{M}\longrightarrow \mathcal{M}'$ is a formal embedding, we have 
 $$\left(\varphi\circ F\circ \psi\right)\left(w_{1},w_{2},\dots,w_{N};z_{1},z_{2},\dots,z_{N}\right)=\left(w_{1},w_{2},\dots,w_{N},0,\dots,0; z_{1},z_{2},\dots,z_{N},0,\dots,0\right),$$
up to compositions with suitable automorphisms  $\varphi$ and $\psi$ of the Models from  (\ref{emo}).
\et 
 
It is an interesting result in the light of the author's previous work\cite{val8}. Its proof is just an application of the methods developed in \cite{val8}. Again, the Fischer Decomposition\cite{sh} is used in order to define Spaces of Fischer Normalizations\cite{sh}.  It simplifies the local defining equations by imposing convenient normalization conditions   with respect to Generalized Fischer Decompositions\cite{sh} similarly as in \cite{val8}. Then,   compositions with suitable automorphisms of the models from (\ref{emo}) gives just a class of equivalence (\ref{emo}). 

 More generally, we consider the following generalized situation motivated by the previous result. Let $\left(\tilde{w},\tilde{z}\right)=\left(\tilde{w}_{1},\tilde{w}_{2},\dots,\tilde{w}_{\tilde{N}};\tilde{z}_{1},\tilde{z}_{2},\dots,\tilde{z}_{\tilde{N}}\right)$ be the coordinates in $\mathbb{C}^{2\tilde{N}}$. Respectively, let $\left(\tilde{w'},\tilde{z}'\right)=\left({\tilde{w'}}_{1},{\tilde{w'}}_{2},\dots,{\tilde{w'}}_{\tilde{N'}};{\tilde{z'}}_{1},{\tilde{z'}}_{2},\dots,{\tilde{z'}}_{\tilde{N'}}\right)$ be the coordinates in $\mathbb{C}^{2\tilde{N'}}$.  
 
 We introduce the following class of Real-Formal Submanifolds in Complex Spaces
  \begin{equation}\mathbb{C}^{2\tilde{N} }\supset M_{\tilde{N}}:\hspace{0.1 cm}\left\{\begin{split}&\tilde{w}_{1}=\tilde{z}_{1}\overline{\tilde{z}}_{1}+\tilde{\lambda}_{1}\left(\tilde{z}_{1}^{2}+\overline{\tilde{z}}_{1}^{2}\right)+\mbox{O}(3) ,\\& \tilde{w}_{2}=\tilde{z}_{2}\overline{\tilde{z}}_{2}+\tilde{\lambda}_{2}\left(\tilde{z}_{2}^{2}+\overline{\tilde{z}}_{2}^{2}\right)+\mbox{O}(3) ,\\&\quad\quad\vdots \\& \tilde{w}_{\tilde{N}}=\tilde{z}_{\tilde{N}}\overline{\tilde{z}}_{\tilde{N}}+\tilde{\lambda}_{\tilde{N}}\left(\tilde{z}_{\tilde{N}}^{2}+\overline{\tilde{z}}_{\tilde{N}}^{2}\right)+\mbox{O}(3) ,\end{split}\right.\quad  \mathbb{C}^{2\tilde{N'} }\supset M_{\tilde{N'}}:\hspace{0.1 cm}\left\{\begin{split}&{\tilde{w'}}_{1}={\tilde{z'}}_{1}\overline{{\tilde{z'}}}_{1}+\tilde{\lambda'}_{1}\left( {\tilde{z}}_{1}^{2}+\overline{\tilde{z'}}_{1}^{2}\right)+\mbox{O}(3) ,\\& {\tilde{w}'}_{2}=\tilde{z'}_{2}\overline{\tilde{z'}}_{2}+\tilde{\lambda'}_{2}\left(\tilde{z'}_{2}^{2}+\overline{\tilde{z'}}_{2}^{2}\right)+\mbox{O}(3) ,\\&\quad\quad\vdots \\& {\tilde{w'}}_{\tilde{N'}}=\tilde{z'}_{\tilde{N'}}\overline{\tilde{z'}}_{\tilde{N'}}+\tilde{\lambda'}_{\tilde{N'}}\left(\tilde{z'}_{\tilde{N'}}^{2}+\overline{\tilde{z'}}_{\tilde{N'}}^{2}\right)+\mbox{O}(3) ,\end{split}\right. \label{emo1}\end{equation}
 where it is assumed the following 
\begin{equation}
\tilde{\lambda}_{1}, \tilde{\lambda}_{2}, \dots,\tilde{\lambda}_{\tilde{N}} \in \left(0,\frac{1}{2}\right),\quad \tilde{\lambda'}_{1}, \tilde{\lambda'}_{2}, \dots,\tilde{\lambda'}_{\tilde{N'}} \in \left(0,\frac{1}{2}\right). \label{lambda1}
 \end{equation}
 
We say that $M_{\tilde{N}}$ is (formally holomorphically) embendable in $\tilde{M}_{N}$ when there exist   formal (holomorphic) embeddings of $M_{\tilde{N}}$ into $\mathcal{M}_{N}$, writing as follows
$$M_{\tilde{N}}\mapsto\mathcal{M}_{N}.$$

 Clearly, this definition is locally considered.  Of course, it remains to study the uniqueness of such embenddings.  We will show  that two such embeddings are  always unique up to a composition with an suitable automorphism of the corresponding model. Then, we can establish similar results  towards to Webster\cite{we1} and Ebenfelt-Huang-Zaitsev\cite{ebhuza} in order to establish  another analogue   of Moser's Theorem\cite{mo}.  
 
 We obtain the following result
\bt\label{tB} Let $\mathcal{M}_{N}\subset\mathbb{C}^{2N}$ and $\mathcal{M}_{N'}'\subset\mathbb{C}^{2N'}$ as in (\ref{emo}), and respectively let  $M_{\tilde{N}}\subset\mathbb{C}^{2\tilde{N}}$ and $M_{\tilde{N'}}\subset\mathbb{C}^{2\tilde{N'}}$ as in (\ref{emo1}). Moreover, we assume that the following holds 

 \begin{equation}\begin{array}[c]{ccc}
 \mathbb{C}^{2N }\supset\mathcal{M}_{N}& \ni\left( w_{1},w_{2},\dots,w_{N},z_{1},z_{2},\dots,z_{N}\right) \longmapsto\left(w_{1},w_{2},\dots,w_{N},0,\dots,0,z_{1},z_{2},\dots,z_{N},0,\dots,0 \right)\in & \mathcal{M}_{N'} \subset \mathbb{C}^{2 N'}\\
 \uparrow\scriptstyle{ } & & \uparrow\scriptstyle{}\\ \mathbb{C}^{2\tilde{N} }
 \supset M_{\tilde{N}} & & M_{\tilde{N'}}\subset \mathbb{C}^{2\tilde{N'} }
 \end{array}. \label{diag1}\end{equation} 

  Then   $M_{\tilde{N}}$ is embedded   into $M_{\tilde{N'}}$ according to the standard linear embedding. Moreover, any formal     embedding  between $M_{\tilde{N}}$ and $M_{\tilde{N'}}$  is equivalent the standard 
  linear embedding. \et 

In particular, according to the last diagram, we obtain 
\begin{equation}
\lambda_{1}=\lambda'_{1}=\tilde{\lambda}_{1}=\tilde{\lambda'}_{1}, \lambda_{2}=\lambda'_{2}=\tilde{\lambda}_{2}=\tilde{\lambda'}_{2},\dots, \lambda_{\tilde{N}}=\lambda'_{\tilde{N}}=\tilde{\lambda}_{\tilde{N}}=\tilde{\lambda'}_{\tilde{N}}.
\end{equation}

\section{Acknowledgements} 

 The Financial Thanks are organized like the relevance of each Funding, emphasizing that  I did not receive the first monthly scholarship   at Federal University of Minas Gerais.  Special Thanks to Science Foundation of Ireland for Generous  Funding during my doctoral period,  and also Generous Thanks  to Capes for Funding  at   Federal University of Santa Catarina, Brazil. It has been a great period for me and I   remember with pleasure the period spent in Florianopolis. I thank also my supervisor (Prof. Dmitri Zaitsev) for many long conversations, regardless that I was a more or less good student. Hopefully, this paper confirms work started then at least partially.

 Regardless that it does not appear written on the published version of \cite{bu1}, I make clear that the reference \cite{bu1} was fully supported by  Science Foundation Ireland Grant 06/RFP/MAT 018, from where I received a scholarship in the first two doctoral years in Trinity College Dublin. It has been an important period for me and I want to defend the Funding existent at the doctoral beginning. 

\section{Ingredients and Computations}
We approach this situation according to  indications from the author's professor\cite{za} recalling   approaches from \cite{bu1},\cite{bu2},\cite{bu3},\cite{bu4},\cite{bu5},\cite{val8}. In particular, we study how a formal (holomorphic)  embedding occurs in the local defining equations (\ref{model2}). We proceed as follows:
\subsection{Settings}Throughout this paper, let  $(w;z)=\left(w_{1},w_{2},\dots,w_{N};z_{1},z_{2},\dots,z_{N}\right)$  be the coordinates in $\mathbb{C}^{2N}$. Clearly, it assumed  $N\in\mathbb{N}^{\star}$.   For the beginning, we shall study a much more general situation:

Let the  permutations $\sigma,\tau\in \mathcal{S}_{N}$. We introduce  the following   model: 
\begin{equation}\mathcal{M}_{N}\left(\sigma,\tau\right): \left\{\begin{split}&w_{1}=z_{1}\overline{z}_{\tau\left(1\right)}+\lambda_{1}\left(z_{1}z_{\sigma\left(1\right)}+\overline{z}_{1}\overline{z}_{\sigma\left(1\right)}\right),\\& w_{2}=z_{2}\overline{z}_{\tau\left(2\right)}+\lambda_{2}\left(z_{2}z_{\sigma \left(2\right)}+\overline{z}_{2}\overline{z}_{\sigma\left(2\right)}\right),\\&\quad\quad\vdots \\& w_{N}=z_{N}\overline{z}_{\tau\left(N\right)}+\lambda_{N}\left(z_{N}z_{\sigma\left(N\right)}+\overline{z}_{N}\overline{z}_{\sigma\left(N\right)}\right).\end{split}\right.\label{model1}  \end{equation}

Let the  permutations $\sigma',\tau'\in \mathcal{S}_{N'}$. We introduce  the following   model: 
\begin{equation}\quad\quad\quad\quad\quad\mathcal{M}_{N'}\left(\sigma',\tau'\right): \left\{\begin{split}&{w'}_{1}={z'}_{1}\overline{{z'}}_{\tau'\left(1\right)}+{\lambda'}_{1}\left({z'}_{1}{z'}_{\sigma'\left(1\right)}+\overline{{z'}}_{1}\overline{z}_{\sigma'\left(1\right)}\right),\\& {w'}_{2}={z'}_{2}\overline{{z'}}_{\tau'\left(2\right)}+{\lambda'}_{2}\left({z'}_{2}{z'}_{\sigma' \left(2\right)}+\overline{{z'}}_{2}\overline{{z'}}_{\sigma'\left(2\right)}\right),\\&\quad\quad\vdots \\& {w'}_{N'}={z'}_{N'}\overline{{z'}}_{\tau'\left(N'\right)}+{\lambda'}_{N'}\left({z'}_{N'}{z'}_{\sigma'\left(N'\right)}+\overline{{z'}}_{N'}\overline{{z'}}_{\sigma\left(N'\right)}\right),\end{split}\right.\label{model2}  \end{equation}
where $\left({w'};{z'}\right)=\left({w'}_{1},{w'}_{2},\dots,{w'}_{N'};{z'}_{1},{z'}_{2},\dots,{z'}_{N'}\right)$  are the coordinates in $\mathbb{C}^{2N'}$. Clearly, it assumed as well $N'\in\mathbb{N}^{\star}$.

\subsection{Formal (Holomorphic) Embeddings} Let $M\subset\mathbb{C}^{2N}$ be a Real-Formal
Submanifold defined near  $p=0$ as follows
\begin{equation}\begin{split} 
w_{l}=z_{l}\overline{z}_{\tau\left(l\right)}+\lambda_{l}\left(z_{l}z_{\sigma\left(l\right)}+\overline{z}_{l}\overline{z}_{\sigma\left(l\right)}\right) +\displaystyle\sum
_{k\geq 3}\varphi_{k}^{\left(l\right)}(z,\overline{z}),\quad\mbox{for all $l=1,\dots,N\leq N'$,}    \end{split} \label{E1}
\end{equation}
where   $\varphi_{k}^{\left(l\right)}(z,\overline{z})$ is a  polynomial of    degree $k$ in $(z,\overline{z})$,  for all  $k\geq 3$ and $l=1,\dots,N\leq N'$.
 
  We consider also another Real-Formal
Submanifold $M'\subset\mathbb{C}^{2N'}$ defined near $p=0$ as follows
\begin{equation}w'_{l'}={z'}_{l'}\overline{{z'}}_{\tau'\left(l'\right)}+{\lambda'}_{l'}\left({z'}_{l'}{z'}_{\sigma'\left(l'\right)}+\overline{{z'}}_{l'}\overline{z}_{\sigma'\left(l'\right)}\right)  +\displaystyle\sum _{k\geq
3}{\varphi'}_{k}^{\left(l'\right)}\left(z',\overline{z'}\right),\quad\mbox{for all $l'=1,\dots,N'$,}\label{E2}
\end{equation}
where ${\varphi'}_{k}^{\left(l'\right)}\left(z',\overline{z'}\right)$ is a
 polynomial of   degree $k$ in $\left(z',\overline{z'}\right)$,  for all  $k\geq 3$ and $l'=1,\dots,N'$.  

We move forward by using the following notations 
\begin{equation}   
Q(z,\overline{z})=\left(q_{1}(z,\overline{z}),q_{2}(z,\overline{z}),\dots,q_{N}(z,\overline{z})\right),  \quad  \mbox{for:} \quad   q_{l}(z,\overline{z})=z_{l}\overline{z}_{\tau\left(l\right)}+\lambda_{l}\left(z_{l}z_{\sigma\left(l\right)}+\overline{z}_{l}\overline{z}_{\sigma\left(l\right)}\right),\label{notation1}
\end{equation}
for all $l=1,\dots,N$, and respectively by using the following notations 
\begin{equation} Q'\left(z',\overline{z'}\right)=\left({q'}_{1}\left(z',\overline{z'}\right),{q'}_{2}\left(z',\overline{z'}\right),\dots,{q'}_{N'}\left(z',\overline{z'}\right)\right),\quad   \mbox{for:} \quad   {q'}_{l'}\left(z',\overline{z'}\right)={z'}_{l'}\overline{{z'}}_{\tau'\left(l'\right)}+{\lambda'}_{l'}\left({z'}_{l'}{z'}_{\sigma'\left(l'\right)}+\overline{{z'}}_{l'}\overline{{z'}}_{\sigma'\left(l'\right)}\right),     \label{notation2}
\end{equation}
for all $l'=1,\dots,N'$.

Let  a formal (holomorphic) embedding of  $M$ into $M'$, denoted as follows
\begin{equation} \left(F(z,w);G(z,w)\right) =\left(F_{1}(z,w),F_{2}(z,w),\dots,F_{N}(z,w),\dots,F_{N'}(z,w) ;G_{1}(z,w),G_{2}(z,w),\dots
,G_{N}(z,w),\dots
,G_{N'}(z,w)\right). \label{gigi}\end{equation} 

Then, for  $w=\left(w_{1},w_{2},\dots,w_{N}\right)$   defined by (\ref{E1}), we have  
\begin{equation} G_{\left(l'\right)}(z,w)=q_{l'}\left(
F(z,w),F(z,w)\right)+\displaystyle\sum _{k\geq
3}{\varphi'}_{k}^{\left(l'\right)}\left(F(z,w),\overline{F(z,w)}\right),\quad\mbox{for all $l'=1,\dots,N'$.} \label{ec}
\end{equation} 
 
In order to understand better the interactions of terms in (\ref{ec}), we write this embedding  (\ref{gigi}) as follows
\begin{equation}\begin{split} \left( F(z,w) ;G(z,w)\right)&  
= \left(\displaystyle\sum_{m;n_{1}, \dots,n_{N}\in\mathbb{N} }F_{m;n_{1}, \dots,n_{N}}^{\left(1\right)}(z)w_{1}^{n_{1}}\dots w_{N}^{n_{N}} ,\dots,\displaystyle\sum_{m;n_{1}, \dots,n_{N}\in\mathbb{N} }F_{m;n_{1}, \dots,n_{N}}^{\left(N\right)}(z)w_{1}^{n_{1}}\dots w_{N}^{n_{N}}, \right.\\&\left.\quad\quad\dots, \displaystyle\sum_{m;n_{1}, \dots,n_{N}\in\mathbb{N} }F_{m;n_{1}, \dots,n_{N}}^{\left(N'\right)}(z)w_{1}^{n_{1}}\dots w_{N}^{n_{N}};\displaystyle\sum_{m;n_{1}, \dots,n_{N}\in\mathbb{N} }G_{m;n_{1}, \dots,n_{N}}^{\left(1\right)}(z)w_{1}^{n_{1}}\dots w_{N}^{n_{N}},\right.\\& \left.\quad\quad\dots,\displaystyle\sum_{m;n_{1}, \dots,n_{N}\in\mathbb{N} }G_{m;n_{1}, \dots,n_{N}}^{\left(N\right)}(z)w_{1}^{n_{1}}\dots w_{N}^{n_{N}} ,\dots,\displaystyle\sum_{m;n_{1}, \dots,n_{N}\in\mathbb{N} }G_{m;n_{1}, \dots,n_{N}}^{\left(N'\right)}(z)w_{1}^{n_{1}}\dots w_{N}^{n_{N}}\right), \end{split}
\label{map} \end{equation}
dealing with the following   homogeneous polynomials of degree
$m$ in $z$
$$G_{m;n_{1}, \dots,n_{N}}^{\left(1\right)}(z),\dots,G_{m;n_{1}, \dots,n_{N}}^{\left(N\right)}(z),\dots,G_{m;n_{1}, \dots,n_{N}}^{\left(N'\right)}(z) ;\quad\quad F_{m;n_{1}, \dots,n_{N}}^{\left(1\right)}(z),\dots,F_{m;n_{1}, \dots,n_{N}}^{\left(N\right)}(z) ,\dots,F_{m;n_{1}, \dots,n_{N}}^{\left(N'\right)}(z),$$  for all  $n_{1}, \dots,n_{N}\in\mathbb{N}$, concluding by (\ref{ec}) and (\ref{map}) the following    $$ G_{0;0,\dots,0}(z)=0,\quad F_{0;0,\dots,0}(z)=0,$$
because the embedding (\ref{gigi}) sends $0\in\mathbb{C}^{2N}$ into $0\in\mathbb{C}^{2N'}$. 
  
 \subsection{Considerations} In order to   understand the appearing interactions of terms in the local defining equations (\ref{E2}) after replacements from (\ref{gigi}), it is important to have a better understanding of the models (\ref{model1}) and \ref{model2}). Thus, let's assume  the following
\begin{equation}\begin{split}&q_{l}(z,\overline{z})=\overline{q_{l}(z,\overline{z})},\quad\mbox{for all $l=1,\dots,N_{0}$,}\\&q_{l}(z,\overline{z})\neq \overline{q_{l}(z,\overline{z})},\quad\mbox{for all $l=N_{0}+1,\dots,N$,}\end{split}\label{F1}
\end{equation}
where $ N_{0}\in  1,\dots, N$.

Similarly as previously, let's assume  the following
\begin{equation}\begin{split}&{q'}_{l'}\left(z',\overline{z'}\right)=\overline{{q'}_{l'}\left(z',\overline{z'}\right)},\quad\hspace{0.1 cm}\mbox{for all $l'=1,\dots,N'_{0}$,}\\&{q'}_{l'}\left(z',\overline{z'}\right)\neq \overline{{q'}_{l'}\left(z',\overline{z'}\right)},\quad\hspace{0.1 cm}\mbox{for all $l'={N'}_{0}+1,\dots,N'$,}\end{split}\label{F2}
\end{equation}
where $ {N'}_{0}\in 1,\dots, N'$.

Then, (\ref{F1}) and (\ref{F2})  yield  that 
\begin{equation}\begin{split}&\tau\left(k\right)=k,\hspace{0.1 cm} \mbox{for all $k=1,\dots,N_{0}$},\quad\quad\quad\tau'\left(k'\right)=k',\hspace{0.15 cm} \mbox{for all $k'=1,\dots,{N'}_{0}$},\\& \tau\left(k\right)\neq k,\hspace{0.1 cm} \mbox{for all $k\not\in 1,\dots,N_{0}$},\quad\quad\quad\tau'\left(k'\right)\neq k',\hspace{0.2 cm} \mbox{for all $k'\not\in 1,\dots,{N'}_{0}$.} \end{split}\label{Fse}
\end{equation}

\subsection{General Equations} By doing  appropriate replacements of (\ref{notation1}),(\ref{notation2}),(\ref{gigi}) in the local defining equations (\ref{E2}),  we  obtain  
\begin{equation}\left.\begin{split}\displaystyle &  \quad\quad\quad\sum
_{m;n_{1}, \dots,n_{N}\in\mathbb{N}}G_{m;n_{1}, \dots,n_{N}}^{\left(l'\right)}(z)\left(q_{1}\left(z,\overline{z}\right) +\displaystyle\sum
_{k\geq 3}\varphi_{k}^{\left(1\right)}(z,\overline{z})\right)^{n_{1}}\dots\left(q_{N'}\left(z,\overline{z}\right) +\displaystyle\sum
_{k\geq 3}\varphi_{k}^{\left(N\right)}(z,\overline{z})\right)^{n_{N}}=\\& \hspace{0.1 cm}\quad q'_{l'}\left(\displaystyle\sum _{m;n_{1}, \dots,n_{N}\in\mathbb{N}}
F_{m;n_{1}, \dots,n_{N}}(z)\left(q_{1}\left(z,\overline{z}\right) +\displaystyle\sum
_{k\geq 3}\varphi_{k}^{\left(1\right)}(z,\overline{z})\right)^{n_{1}}\dots\left(q_{N}\left(z,\overline{z}\right) +\displaystyle\sum
_{k\geq 3}\varphi_{k}^{\left(N\right)}(z,\overline{z})\right)^{n_{N}}\right.,\\&\left.\quad\quad\quad\overline{\displaystyle\sum _{m;n_{1}, \dots,n_{N}\in\mathbb{N}}
F_{m;n_{1}, \dots,n_{N}}(z)\left(q_{1}\left(z,\overline{z}\right) +\displaystyle\sum
_{k\geq 3}\varphi_{k}^{\left(1\right)}(z,\overline{z})\right)^{n_{1}}\dots\left(q_{N}\left(z,\overline{z}\right) +\displaystyle\sum
_{k\geq 3}\varphi_{k}^{\left(N\right)}(z,\overline{z})\right)^{n_{N}}}\right)  \\&+\displaystyle\sum
_{k\geq 3}{\varphi'}_{k}^{\left(l'\right)}
\left(\displaystyle\sum _{m;n_{1}, \dots,n_{N}\in\mathbb{N}}F_{m,n}(z)\left(q_{1}\left(z,\overline{z}\right) +\displaystyle\sum
_{k\geq 3}\varphi_{k}^{\left(1\right)}(z,\overline{z})\right)^{n_{1}}\dots\left(q_{N}\left(z,\overline{z}\right) +\displaystyle\sum
_{k\geq 3}\varphi_{k}^{\left(N\right)}(z,\overline{z})\right)^{n_{N}}\right.,\\& \quad\quad\quad\quad\quad\quad\left.\overline{\displaystyle\sum _{m;n_{1}, \dots,n_{N}\in\mathbb{N}}F_{m,n}(z)\left(q_{1}\left(z,\overline{z}\right) +\displaystyle\sum
_{k\geq 3}\varphi_{k}^{\left(1\right)}(z,\overline{z})\right)^{n_{1}}\dots\left(q_{N}\left(z,\overline{z}\right) +\displaystyle\sum
_{k\geq 3}\varphi_{k}^{\left(N\right)}(z,\overline{z})\right)^{n_{N}}}\right),
\end{split}\right.
\label{ecuatiegenerala}\end{equation}
for all $l'=1,\dots,N'$. 

\subsection{Linear Changes of Coordinates}We can  assume  the following  \begin{equation}\left(F_{_{1;0,\dots,0}}^{\left(1\right)}(z),\dots,F_{_{1;0,\dots,0}}^{\left(N\right)}(z),F_{_{1;0,\dots,0}}^{\left(N+1\right)}(z),\dots,F_{_{1;0,\dots,0}}^{\left(N'\right)}(z)\right)=\left(z_{1},\dots,z_{N},0,\dots,0\right),\label{bolivia1}\end{equation}
and respectively, we can assume the following
 \begin{equation} \begin{pmatrix} G^{\left(1\right)}_{0;1,0,\dots,0}(z) & G_{0;0 ,1,\dots,0}^{\left(1\right)}(z)& \dots & G_{0;0,1,\dots,0}^{\left(1\right)}(z)\\ G^{\left(2\right)}_{0;1,0,\dots,0}(z) & G_{0;0,1,\dots,0}^{\left(2\right)}(z)& \dots & G_{0;0,\dots,1}^{\left(2\right)}(z)  \\  \vdots &\vdots &\ddots&\vdots   \\ G^{\left(N\right)}_{0;1,0,\dots,0}(z) & G_{0;0,1,\dots,0}^{\left(N\right)}(z)& \dots & G_{0;0,0,\dots,1}^{\left(N\right)}(z) \\  \vdots &\vdots &\ddots&\vdots    \\  G^{\left(N'\right)}_{0;1,0,\dots,0}(z) & G_{0;0,1,\dots,0}^{\left(N'\right)}(z)& \dots & G_{0;0,\dots,1}^{\left(N'\right)}(z) \end{pmatrix}=\begin{pmatrix} 1 & 0 & \dots & 0    \\ 0 & 1 & \dots & 0   \\  \vdots &\vdots &\ddots&\vdots    \\ 0 & 0 & \dots & 1    \\  \vdots &\vdots &\ddots&\vdots    \\ 0 & 0 & \dots & 0    \end{pmatrix} ,\label{bolivia2}\end{equation} 
 by   composing (\ref{map}) with a   linear holomorphic automorphism of the    models (\ref{model1}) and (\ref{model2}). 
 
 Indeed, if we write as follows 
 \begin{equation}\begin{pmatrix} F_{_{1;0,\dots,0}}^{\left(1\right)}(z)   \\ F_{_{1;0,\dots,0}}^{\left(2\right)}(z) \\  \vdots &    \\ F_{_{1;0,\dots,0}}^{\left(N\right)}(z)  \\  \vdots &    \\ F_{_{1;0,\dots,0}}^{\left(N'\right)}(z)\end{pmatrix} =\begin{pmatrix} \gamma_{11} & \gamma_{12} & \dots & \gamma_{1N}  \\ \gamma_{21} & \gamma_{22} & \dots & \gamma_{2N}   \\  \vdots &\vdots &\ddots&\vdots    \\ \gamma_{N1} & \gamma_{N2} & \dots & \gamma_{NN} \\  \vdots &\vdots &\ddots&\vdots   \\ \gamma_{N'1} & \gamma_{N'2} & \dots & \gamma_{N'N} \end{pmatrix} \begin{pmatrix} z_{1}   \\ z_{2} \\  \vdots &    \\ z_{N}\end{pmatrix},
\label{bolivia3}\end{equation}
where we deal with 
$$\gamma_{11},\dots,\gamma_{1N},\dots,\gamma_{1N'}; \gamma_{21},\dots,\gamma_{2N},\dots,\gamma_{2N'};\dots;\gamma_{N1},\dots,\gamma_{NN},\dots,\gamma_{NN'};\dots; \gamma_{N'1},\dots,\gamma_{N'N},\dots,\gamma_{N'N'}\in\mathbb{R}, 
$$
we obtain an system of  equations  written in its matrix form as follows
\begin{equation}\begin{split}&\quad\quad\quad\quad\quad\quad \begin{pmatrix} G^{\left(1\right)}_{0;1,0,\dots,0}(z) & G_{0;0 ,1,\dots,0}^{\left(1\right)}(z)& \dots & G_{0;0,1,\dots,0}^{\left(1\right)}(z)\\ G^{\left(2\right)}_{0;1,0,\dots,0}(z) & G_{0;0,1,\dots,0}^{\left(2\right)}(z)& \dots & G_{0;0,\dots,1}^{\left(2\right)}(z)  \\  \vdots &\vdots &\ddots&\vdots   \\ G^{\left(N\right)}_{0;1,0,\dots,0}(z) & G_{0;0,1,\dots,0}^{\left(N\right)}(z)& \dots & G_{0;0,0,\dots,1}^{\left(N\right)}(z) \\  \vdots &\vdots &\ddots&\vdots    \\  G^{\left(N'\right)}_{0;1,0,\dots,0}(z) & G_{0;0,1,\dots,0}^{\left(N'\right)}(z)& \dots & G_{0;0,\dots,1}^{\left(N'\right)}(z) \end{pmatrix}\begin{pmatrix} z_{1}\overline{z}_{\tau\left(1\right)}+\lambda_{1}\left(z_{1}z_{\sigma\left(1\right)}+\overline{z}_{1}\overline{z}_{\sigma\left(l\right)}\right)   \\ z_{2}\overline{z}_{\tau\left(2\right)}+\lambda_{2}\left(z_{2}z_{\sigma\left(2\right)}+\overline{z}_{2}\overline{z}_{\sigma\left(2\right)}\right) \\  \vdots &\   \\ z_{N}\overline{z}_{\tau\left(N\right)}+\lambda_{N}\left(z_{N}z_{\sigma\left(N\right)}+\overline{z}_{N}\overline{z}_{\sigma\left(N\right)}\right)\end{pmatrix} =\\&\begin{pmatrix}
\left(\displaystyle\sum_{i=1}^{N}\gamma_{1i}z_{i}\right)\overline{\left(\displaystyle\sum_{i=1}^{N}\gamma_{ \tau'\left(1\right)i}z_{i}\right)} +{\lambda'}_{1}\left(\left(\displaystyle\sum_{i=1}^{N}\gamma_{1i}z_{i}\right)\left(\displaystyle\sum_{i=1}^{N}\gamma_{ \sigma'\left(1\right)i}z_{i}\right) +\overline{\left(\displaystyle\sum_{i=1}^{N}\gamma_{1i}z_{i}\right)} \overline{\left(\displaystyle\sum_{i=1}^{N}\gamma_{ \sigma'\left(1\right)i}z_{i}\right)} \right)   \\ 
 \left(\displaystyle\sum_{i=1}^{N}\gamma_{2i}z_{i}\right)
\overline{\left(\displaystyle\sum_{i=1}^{N}\gamma_{ \tau'\left(2\right)i}z_{i}\right)} +{\lambda'}_{2}\left(\left(\displaystyle\sum_{i=1}^{N}\gamma_{2i}z_{i}\right)
\left(\displaystyle\sum_{i=1}^{N}\gamma_{ \sigma'\left(2\right)i}z_{i}\right)+\overline{ \left(\displaystyle\sum_{i=1}^{N}\gamma_{2i}z_{i}\right)
} \overline{\left(\displaystyle\sum_{i=1}^{N}\gamma_{ \sigma'\left(1\right)i}z_{i}\right)} \right)  \\  \vdots &\   \\ 
\left(\displaystyle\sum_{i=1}^{N}\gamma_{N'i}z_{i}\right)\overline{\left(\displaystyle\sum_{i=1}^{N}\gamma_{\tau'\left(N'\right)i}z_{i}\right)}+{\lambda'}_{N'}\left(\left(\displaystyle\sum_{i=1}^{N}\gamma_{N'i}z_{i}\right)  \left(\displaystyle\sum_{i=1}^{N}\gamma_{ \sigma'\left(N'\right)i}z_{i}\right)+\overline{\left(\displaystyle\sum_{i=1}^{N}\gamma_{N'i}z_{i}\right)} \overline{\left(\displaystyle\sum_{i=1}^{N}\gamma_{ \sigma'\left(N'\right)i}z_{i}\right)} \right)\end{pmatrix},\end{split} \label{sistem1}
\end{equation} 
having a matrix of type ${N'}\times N$ in the above left-hand side.  

Moreover, we can assume above the following
\begin{equation}\det  \begin{pmatrix} G^{\left(1\right)}_{0;1,0,\dots,0}(z) & G_{0;0 ,1,\dots,0}^{\left(1\right)}(z)& \dots & G_{0;0,1,\dots,0}^{\left(1\right)}(z)\\ G^{\left(2\right)}_{0;1,0,\dots,0}(z) & G_{0;0,1,\dots,0}^{\left(2\right)}(z)& \dots & G_{0;0,\dots,1}^{\left(2\right)}(z)  \\  \vdots &\vdots &\ddots&\vdots   \\ G^{\left(N\right)}_{0;1,0,\dots,0}(z) & G_{0;0,1,\dots,0}^{\left(N\right)}(z)& \dots & G_{0;0,0,\dots,1}^{\left(N\right)}(z)   \end{pmatrix} \neq 0,\label{gig4}
\end{equation}
because otherwise we can linearly change the coordinates  $\left(w_{1},w_{2},\dots,w_{N}\right)$ in order to assume that (\ref{gig4}) holds, given that the entries of the above left-hand side matrix  are just numbers.
 
Now, extracting the terms of degree $2$ from (\ref{ecuatiegenerala}) or identifying the coefficients in (\ref{sistem1}), we obtain
\begin{equation} G^{\left(l'\right)}_{0;1,0,\dots,0}(z)z_{1}\overline{z}_{\tau\left(1\right)} + G_{0;0 ,1,\dots,0}^{\left(l'\right)}(z)z_{2}\overline{z}_{\tau\left(2\right)}+ \dots + G_{0;0,1,\dots,0}^{\left(l'\right)}(z)z_{N}\overline{z}_{\tau\left(N\right)} =\left(\displaystyle\sum_{i=1}^{N}\gamma_{l'i}z_{i}\right)\overline{\left(\displaystyle\sum_{i=1}^{N}\gamma_{ \tau'\left(l'\right)i}z_{i}\right)}, \label{iri1}
\end{equation}
for all $l'=1,\dots,N'$.

Analogously as above, we obtain
\begin{equation} G^{\left(l'\right)}_{0;1,0,\dots,0}(z)z_{1}z_{\sigma\left(1\right)} + G_{0;0 ,1,\dots,0}^{\left(l'\right)}(z)z_{2}z_{\sigma\left(2\right)}+ \dots + G_{0;0,1,\dots,0}^{\left(l'\right)}(z)z_{N}z_{\sigma\left(N\right)} =\frac{{\lambda'}_{l'}}{\lambda_{l'}}\left(\displaystyle\sum_{i=1}^{N}\gamma_{l'i}z_{i}\right)\left(\displaystyle\sum_{i=1}^{N}\gamma_{ \sigma'\left(l'\right)i}z_{i}\right),\label{iri2}
\end{equation}
for all $l'=1,\dots,N'$.

Now, (\ref{iri1}) yields that
 \begin{equation}\gamma_{l'i}\overline{\gamma_{\tau'\left(l'\right)j}}=\left\{ \begin{matrix} G^{\left(l'\right)}_{0;0, \dots,1,\dots, 0}(z),&\quad \mbox{ for all $i,j=1,\dots,N$ with $j= \tau\left(i\right) $,}  \\ \quad\quad\quad\quad\quad\quad 0,&\quad\mbox{ for all $i,j=1,\dots,N$ with $j\neq   \tau\left(i\right) $,}
  \end{matrix} \right.\label{gig1}\end{equation}
for all $l'=1,\dots,N'$.
  
On the other hand, (\ref{iri2}) yields that
 \begin{equation}\frac{{\lambda'}_{l'}}{\lambda_{l'}}\gamma_{l'i} \gamma_{\sigma'\left(l'\right)j} =\left\{ \begin{matrix} G^{\left(l'\right)}_{0;0, \dots,1,\dots, 0}(z),&\quad \mbox{ for all $i,j=1,\dots,N$ with $j= \sigma\left(i\right) $,}  \\ \quad\quad\quad\quad\quad\quad 0,&\quad\mbox{ for all $i,j=1,\dots,N$ with $j\neq   \sigma\left(i\right) $,}
  \end{matrix} \right.\label{gig2}\end{equation}
for all $l'=1,\dots,N'$.

Let $l'\in 1,\dots, N'$ and $i_{1},i_{2}\in 1,\dots,N$ such that $i_{2} \neq \tau\left(i_{1}\right)$, satisfying that
\begin{equation}\gamma_{l'i_{1}}\neq 0,\quad \gamma_{\tau'\left(l'\right) i_{2} }\neq 0.\label{gig3}\end{equation} 

 Then, according to the above observations (\ref{gig1}), we have
\begin{equation}\gamma_{l'i_{1}} \gamma_{\tau'\left(l'\right) i_{2} }= 0,\label{gig3A}\end{equation}
which gives a contraction, resulting therefore  that $\tau\left(i_{1}\right)=i_{2}$, if (\ref{gig3}) holds.  

It is obvious the  existence  of such numbers as $i_{1}$, because (\ref{gigi}) is an embedding and therefore (\ref{gig4}) holds. Moreover,  the following conclusion becomes by (\ref{gig1}) clear: there exists a unique $i_{l}\in 1,\dots,N$ such that the following holds
\begin{equation}0\neq \gamma_{l'i}\overline{\gamma_{\tau'\left(l'\right)\tau(i)}} =G^{\left(1\right)}_{0;0, \dots,1,\dots, 0}(z)\Longleftrightarrow i=i_{l},\quad\mbox{for all $l=1,\dots,N$.} \label{gig6}
\end{equation}

On the other hand from (\ref{gig2}), we have 
\begin{equation}\frac{{\lambda'}_{l}}{\lambda_{l}}\gamma_{li} \gamma_{\sigma'\left(l\right)\sigma(i)}=G^{\left(1\right)}_{0;0, \dots,1,\dots, 0}(z),\quad\mbox{for all $l=1,\dots,N$.}\label{gig7}
\end{equation}

Now, combining (\ref{gig4}),(\ref{iri1}),(\ref{iri2}),(\ref{gig1}),(\ref{gig2}),(\ref{gig3}),(\ref{gig6}),(\ref{gig7}), we obtain 
\begin{equation*}\det \left(\left(\gamma_{li}\right)_{1\leq l\leq N\atop1\leq i\leq N }\right)\neq 0,
\end{equation*}
because (\ref{gig51}) holds, according to the following equalities 
\begin{equation}\begin{split}&\begin{pmatrix} G^{\left(1\right)}_{0;1,0,\dots,0}(z) & G_{0;0 ,1,\dots,0}^{\left(1\right)}(z)& \dots & G_{0;0,1,\dots,0}^{\left(1\right)}(z)\\ G^{\left(2\right)}_{0;1,0,\dots,0}(z) & G_{0;0,1,\dots,0}^{\left(2\right)}(z)& \dots & G_{0;0,\dots,1}^{\left(2\right)}(z)  \\  \vdots &\vdots &\ddots&\vdots   \\ G^{\left(N\right)}_{0;1,0,\dots,0}(z) & G_{0;0,1,\dots,0}^{\left(N\right)}(z)& \dots & G_{0;0,0,\dots,1}^{\left(N\right)}(z)   \end{pmatrix}=\begin{pmatrix} \gamma_{11} & \gamma_{12} & \dots & \gamma_{1N}    \\ \gamma_{21} & \gamma_{22} & \dots & \gamma_{2N}   \\  \vdots &\vdots &\ddots&\vdots    \\ \gamma_{N1} & \gamma_{N2} & \dots & \gamma_{NN}      \end{pmatrix}\begin{pmatrix} \overline{\gamma}_{\tau'(1)1} & \overline{\gamma}_{\tau'(2)1} & \dots & \overline{\gamma}_{\tau'\left(N\right)1}    \\ \overline{\gamma}_{\tau'(1)2} & \overline{\gamma}_{\tau'(2)2} & \dots & \overline{\gamma}_{\tau'\left(N\right)2}   \\  \vdots &\vdots &\ddots&\vdots    \\ \overline{\gamma}_{\tau'\left(1\right)N} & \overline{\gamma}_{\tau'\left(2\right)N} & \dots & \overline{\gamma}_{\tau'\left(N\right)N}      \end{pmatrix}\\&\quad\quad\quad\quad\quad\quad\quad\quad\quad\quad\quad\quad=\begin{pmatrix} \frac{{\lambda'}_{1}}{\lambda_{1}} & 0 & \dots & 0    \\ 0 & \frac{{\lambda'}_{2}}{\lambda_{2}} & \dots & 0   \\  \vdots &\vdots &\ddots&\vdots    \\ 0 & 0 & \dots & \frac{{\lambda'}_{N}}{\lambda_{N}}      \end{pmatrix}\begin{pmatrix} \gamma_{11} & \gamma_{12} & \dots & \gamma_{1N}    \\ \gamma_{21} & \gamma_{22} & \dots & \gamma_{2N}   \\  \vdots &\vdots &\ddots&\vdots    \\ \gamma_{N1} & \gamma_{N2} & \dots & \gamma_{NN}      \end{pmatrix}\begin{pmatrix} \gamma_{\sigma'(1)1} & \gamma_{\sigma'(2)1} & \dots & \gamma_{\sigma'\left(N\right)1}    \\ \gamma_{\sigma'\left(2\right)1} & \gamma_{\sigma'(2)2} & \dots & \gamma_{\sigma'\left(N\right)2}   \\  \vdots &\vdots &\ddots&\vdots    \\ \gamma_{\sigma'\left(1\right)N} & \gamma_{\sigma'\left(2\right)N} & \dots & \gamma_{\sigma'\left(N\right)N}      \end{pmatrix},\end{split} \label{gig51}
\end{equation}
keeping in mind that just an entry does not vanish on each row of the above left-hand matrices, and resulting also that a similar property is satisfied by each of the matrices from the above right-hand matrices. In particular, we also obtain 
\begin{equation}\lambda_{1}={\lambda'}_{1},\lambda_{2}={\lambda'}_{2},\dots,\lambda_{N}={\lambda'}_{N},\quad\tau\equiv \sigma,\quad 0\neq \gamma_{li}\overline{\gamma_{\tau'\left(l\right)\tau(i)}}=\gamma_{li}\gamma_{\tau'\left(l\right)\tau(i)},\hspace{0.1 cm}\mbox{for all $l,i=1,\dots,N$},\label{opa1}
\end{equation}
and therefore, we have
\begin{equation}\gamma_{li}=\overline{\gamma_{li}},\hspace{0.1 cm}\mbox{for all $l,i=1,\dots,N$.}\label{opa2}
\end{equation} 

Moreover, it follows that 
$$ N_{0}\leq {N'}_{0}. $$
 
Now,   we consider the following change of coordinates
\begin{equation}\begin{pmatrix} {w'}_{1} \\ {w'}_{2}\\ \vdots \\ {w'}_{N}
\end{pmatrix}=\begin{pmatrix} G^{\left(1\right)}_{0;1,0,\dots,0}(z) & G_{0;0 ,1,\dots,0}^{\left(1\right)}(z)& \dots & G_{0;0,1,\dots,0}^{\left(1\right)}(z)\\ G^{\left(2\right)}_{0;1,0,\dots,0}(z) & G_{0;0,1,\dots,0}^{\left(2\right)}(z)& \dots & G_{0;0,\dots,1}^{\left(2\right)}(z)  \\  \vdots &\vdots &\ddots&\vdots   \\ G^{\left(N\right)}_{0;1,0,\dots,0}(z) & G_{0;0,1,\dots,0}^{\left(N\right)}(z)& \dots & G_{0;0,0,\dots,1}^{\left(N\right)}(z)   \end{pmatrix}^{-1}\begin{pmatrix} w_{1} \\ w_{2}\\ \vdots \\ w_{N}
\end{pmatrix},\quad \begin{pmatrix} {z'}_{1} \\ {z'}_{2}\\ \vdots \\ {z'}_{N}
\end{pmatrix}=\begin{pmatrix} \gamma_{11} & \gamma_{12} & \dots & \gamma_{1N}    \\ \gamma_{21} & \gamma_{22} & \dots & \gamma_{2N}   \\  \vdots &\vdots &\ddots&\vdots    \\ \gamma_{N1} & \gamma_{N2} & \dots & \gamma_{NN}      \end{pmatrix}^{-1}\begin{pmatrix} z_{1} \\ z_{2}\\ \vdots \\ z_{N}
\end{pmatrix}, \label{opa3} 
\end{equation}
which clearly preserves the model from (\ref{model2}) in the light of (\ref{gig1}),(\ref{gig2}),(\ref{opa1}) and (\ref{opa2}).

Therefore,   changing linearly the coordinates in $\left(w_{1},w_{2},\dots,w_{N}\right)$ and   $\left(z_{1},z_{2},\dots,z_{N}\right)$ according to (\ref{opa3}), we can   assume that 
\begin{equation}\begin{split}&\begin{pmatrix} G^{\left(1\right)}_{0;1,0,\dots,0}(z) & G_{0;0 ,1,\dots,0}^{\left(1\right)}(z)& \dots & G_{0;0,1,\dots,0}^{\left(1\right)}(z)\\ G^{\left(2\right)}_{0;1,0,\dots,0}(z) & G_{0;0,1,\dots,0}^{\left(2\right)}(z)& \dots & G_{0;0,\dots,1}^{\left(2\right)}(z)  \\  \vdots &\vdots &\ddots&\vdots   \\ G^{\left(N\right)}_{0;1,0,\dots,0}(z) & G_{0;0,1,\dots,0}^{\left(N\right)}(z)& \dots & G_{0;0,0,\dots,1}^{\left(N\right)}(z)   \end{pmatrix}=\begin{pmatrix} 1 & 0 & \dots & 0    \\ 0 & 1 & \dots & 0   \\  \vdots &\vdots &\ddots&\vdots    \\ 0 & 0 & \dots & 1      \end{pmatrix},\\&\begin{pmatrix} F^{\left(1\right)}_{0;1,0,\dots,0}(z) & F_{0;0 ,1,\dots,0}^{\left(1\right)}(z)& \dots & F_{0;0,1,\dots,0}^{\left(1\right)}(z)\\ F^{\left(2\right)}_{0;1,0,\dots,0}(z) & F_{0;0,1,\dots,0}^{\left(2\right)}(z)& \dots & F_{0;0,\dots,1}^{\left(2\right)}(z)  \\  \vdots &\vdots &\ddots&\vdots   \\ F^{\left(N\right)}_{0;1,0,\dots,0}(z) & F_{0;0,1,\dots,0}^{\left(N\right)}(z)& \dots & F_{0;0,0,\dots,1}^{\left(N\right)}(z)   \end{pmatrix}=\begin{pmatrix} 1 & 0 & \dots & 0    \\ 0 & 1 & \dots & 0   \\  \vdots &\vdots &\ddots&\vdots    \\ 0 & 0 & \dots & 1      \end{pmatrix}. \end{split} \label{gig5}
\end{equation}
 
We move forward. We have
\begin{equation}\displaystyle\sum_{i=1}^{N}G^{\left(l\right)}_{0;0,\dots,1,0\dots,0}(z)w_{i}=\left(\displaystyle\sum_{i=1}^{N}\gamma_{li}z_{i}\right)\overline{\left(\displaystyle\sum_{i=1}^{N}\gamma_{\tau'\left(l\right)i}z_{i}\right)}+{\lambda'}_{l}\left(\left(\displaystyle\sum_{i=1}^{N}\gamma_{li}z_{i}\right)  \left(\displaystyle\sum_{i=1}^{N}\gamma_{ \sigma'\left(N'\right)i}z_{i}\right)+\overline{\left(\displaystyle\sum_{i=1}^{N}\gamma_{li}z_{i}\right)} \overline{\left(\displaystyle\sum_{i=1}^{N}\gamma_{ \sigma'\left(l\right)i}z_{i}\right)} \right),\label{iri} 
\end{equation}
for all $l=N+1,\dots, N'$.

Now, combining the arguments related to (\ref{gig3}) and (\ref{gig3A}), we obtain
 \begin{equation}\gamma_{l'i}\overline{\gamma_{\tau'\left(l'\right)j}}=\left\{ \begin{matrix} G^{\left(l'\right)}_{0;0, \dots,1,\dots, 0}(z),&\quad \mbox{ for all $i=1,\dots,N$ and $j=N+1,\dots,N'$ with $j= \tau\left(i\right) $,}  \\ \quad\quad\quad\quad\quad\quad 0,&\quad\mbox{ for all $i=1,\dots,N$ and $j=N+1,\dots,N'$ with $j\neq   \tau\left(i\right) $,}
  \end{matrix} \right.\label{gig1A}\end{equation}
and respectively, we obtain 
\begin{equation}\frac{{\lambda'}_{l'}}{\lambda_{l'}}\gamma_{l'i} \gamma_{\sigma'\left(l'\right)j} =\left\{ \begin{matrix} G^{\left(l'\right)}_{0;0, \dots,1,\dots, 0}(z),&\quad \mbox{ for all $i=1,\dots,N$ and $j=N+1,\dots,N'$ with $j= \sigma\left(i\right) $,}  \\ \quad\quad\quad\quad\quad\quad 0,&\quad\mbox{ for all $i=1,\dots,N$ and $j=N+1,\dots,N'$ with $j\neq   \sigma\left(i\right) $,}
  \end{matrix} \right.\label{gig2A}\end{equation}
for all $l'=N+1,\dots,N'$.

Now, we shall pursue several consecutive changes of coordinates in order to normalize the mapping (\ref{gigi}). We proceed as follows: we assume that $l'=N+1$. Moreover, we assume that 
\begin{equation} G^{\left(N+1\right)}_{0;1,0, \dots,0,\dots, 0}(z)=G^{\left(N+1\right)}_{0;0,0, \dots,1,\dots, 0}(z)=\dots=G^{\left(N+1\right)}_{0;0,0, \dots,0,\dots, 1}(z)=0,\label{alfa1}\end{equation}
obtaining that 
\begin{equation}\gamma_{N+1,1}=\gamma_{N+1,2}=\dots=\gamma_{N+1,N}=0.\label{alfa2}\end{equation}
 
 Contrary, if we have for example the following case $$G^{\left(N+1\right)}_{0;1,0, \dots,0,\dots, 0}(z)\neq 0,$$
 then, exactly as previously, we obtain 
  \begin{equation}\lambda_{1}={\lambda'}_{N+1},\quad\tau(1)\equiv \sigma(1),\quad 0\neq \gamma_{N+1,1}\overline{\gamma_{\tau'\left(N+1\right)\tau(1)}}=\gamma_{N+1,1}\gamma_{\sigma'\left(N+1\right)\tau(1)}=G^{\left(N+1\right)}_{0;1,0, \dots,0,\dots, 0}(z),\quad\gamma_{N+1,1}=\overline{\gamma_{N+1,1}}. \label{opa3} 
\end{equation}
In such case we consider a change of coordinates in $(z,w)$ like the following
  \begin{equation}\begin{split}&w'=\left({w'}_{1},\dots,{w'}_{N},\frac{{w'}_{N+1}}{G^{\left(N+1\right)}_{0;1,0, \dots,0,\dots, 0}(z)},{w'}_{N+2},\dots, {w'}_{N'}\right), \\&{z'}=\left( {z'}_{1},\dots,{z'}_{N},\frac{{z'}_{N+1}}{\gamma_{N+1,1}},{z'}_{N+2},\dots, {z'}_{N'}\right) , \label{opa4}\end{split}\end{equation}
which obviously preserves the model (\ref{model2}), because of (\ref{opa3}). 

 Now, in these new coordinates (\ref{opa4}) we can assume that (\ref{alfa1}) and (\ref{alfa2}) hold by the following change  of coordinates
  \begin{equation}\begin{split}&w'=\left({w'}_{1},\dots,{w'}_{N}, {w'}_{N+1}-w_{1},{w'}_{N+2},\dots, {w'}_{N'}\right) \\&{z'}=\left( {z'}_{1},\dots,{z'}_{N}, {z'}_{N+1} ,{z'}_{N+2},\dots, {z'}_{N'}\right) , \end{split}\label{opa5}\end{equation}
which obviously preserves the model (\ref{model2}), because of (\ref{opa3}). 

Therefore, by repeating the previous procedure related to (\ref{opa4}) and (\ref{opa5}), we can make several consecutive linear changes of coordinates in order to normalize the embedding (\ref{gigi}), denoted as in (\ref{map}), as follows
\begin{equation} \begin{split}& F(z,w)  =\left(z_{1}+\mbox{O}(2),z_{2}+\mbox{O}(2),\dots,z_{N}+\mbox{O}(2),\mbox{O}(2),\dots,\mbox{O}(2)\right),\\&  G(z,w)  =\left( w_{1}+\mbox{O}(2),w_{2}+\mbox{O}(2),\dots,w_{N}+\mbox{O}(2),\mbox{O}(2),\dots,\mbox{O}(2)\right). \end{split}\label{gigin}\end{equation}

Next, we continue to study the interactions of the homogeneoys terms in (\ref{ecuatiegenerala}), which are much more complicated than in \cite{bu1}, being similar to the computations from \cite{val8}. Thus, the  computational obstacles are eliminated applying the approach developed by the author\cite{val8}. More precisely, we consider iterative Fischer Decompositions\cite{sh} as follows:

\subsection{Fischer Decompositions\cite{sh}}Recalling the strategy from \cite{bu2}, we define  
\begin{equation} \tr_{l}=\frac{\partial^{2}}{\partial z_{l}\partial \overline{z}_{l}}+\lambda_{l}\left(\frac{\partial^{2}}{\partial z_{l}^{2}}+\frac{\partial^{2}}{\partial \overline{z}_{l}^{2}  }\right),\quad\mbox{for all $l=1,\dots,N$}. \label{tracce}
\end{equation}
 
These differential operators (\ref{tracce})  are just the Fischer Differential Operators associated to  homogeneous polynomials from the right-hand sides in (\ref{model1}),  which may not be  real as it was explained earlier. This fact   influences the further Generalized Fischer Decompositions applied from Shapiro\cite{sh} according the approach\cite{bu1} learned by the author from  Zaitsev\cite{za}. 
We   write  as follows
\begin{equation} P(z)=\displaystyle\sum_{l=1}^{N}A_{l}(z,\overline{z})q_{l}(z,\overline{z}) +P_{0}(z,\overline{z}),\label{V1} \end{equation}
where the following is satisfied
\begin{equation} P_{0}(z,\overline{z})\in   
\displaystyle\bigcap_{l=1}^{N}\ker \tr_{l}  ,\label{V1s} \end{equation}
for any given homogeneous polynomial  $P(z)$.

These resulted homogeneous polynomials,  from (\ref{V1}), are used in order to apply the generalized version of the Fischer Decomposition\cite{sh}   by separating the real parts and the imaginary parts of the local defining equations at each degree level recalling the strategy from \cite{val8}. In particular, we use    the following notations:
\begin{equation} \begin{split}&\hspace{0.03 cm}
z^{I}=z^{i_{1}}\cdot\dots\cdot z^{i_{N}},\quad\mbox{where $\hspace{0.03 cm}I=\left(i_{1},\dots,i_{N}\right)\in\mathbb{N}^{N}$ such that $\hspace{0.1 cm}\left|I\right|=i_{1}+\dots+i_{N}\geq 3$,}\\& z^{J}=z^{j_{1}}\cdot\dots\cdot z^{j_{N}},\quad\mbox{where $J=\left(j_{1},\dots,j_{N}\right)\in\mathbb{N}^{N}$ such that $\left|J\right|=j_{1}+\dots+j_{N}\geq 2$.} \end{split} \label{yy}
\end{equation}

We move forward:
\section{Spaces of Fischer $G$-Decompositions \cite{bu2},\cite{bu3}}We consider by (\ref{model1}) the following Generalized Fischer Decomposition 
\begin{equation}z^{\tilde{I}}=\displaystyle\sum_{l=1}^{N}A_{l}(z,\overline{z}) q_{l}(z,\overline{z})+C(z,\overline{z}),\quad \displaystyle\bigcap_{l=1}^{N} \tr_{l}\left(C(z,\overline{z})\right)=0, \label{op}
\end{equation}
where  $\tilde{I}=\left(\tilde{i}_{1},\tilde{i}_{2},\dots,\tilde{i}_{N}\right)\in\mathbb{N}^{N}$ and $\tilde{i}_{1}+\tilde{i}_{2}+\dots+\tilde{i}_{N}=p$. 
 
We compute these homogeneous polynomials $A_{1}(z,\overline{z}),A_{2}(z,\overline{z}),\dots, A_{N}(z,\overline{z})$ and $C(z,\overline{z})$  straightforwardly  in (\ref{op}),  writing   as follows
\begin{equation} A_{l}(z,\overline{z})=\displaystyle\sum_{I,J\in\mathbb{N}^{N}\atop\left|I\right|+\left|J\right|=p-2} a_{I;J}^{(l)}z^{I}\overline{z}^{J}, \quad C(z,\overline{z}) =\displaystyle\sum_{I,J\in\mathbb{N}^{N}\atop\left|I\right|+\left|J\right|=p} c_{I;J}z^{I}\overline{z}^{J},\quad\mbox{for all $l\in 1,\dots, N$.}
\label{oppv}
\end{equation}

Following \cite{val8}, we make computations differentiating in (\ref{op}), and respectively in (\ref{oppv}). We obtain
 \begin{equation}\begin{split}\quad\quad\quad\lambda_{l} i_{l}\left(i_{l}-1\right)z_{1}^{\tilde{i}_{1}}\dots z_{l}^{\tilde{i}_{l}-2}\dots z_{N}^{\tilde{i}_{N}} =&\displaystyle\sum_{l=1}^{N}\tr_{l}\left( A_{l}(z,\overline{z})\right) q_{l}(z,\overline{z})+\left(z_{l}\left(1+4\lambda_{l}^{2}\right)+4\lambda_{l}\overline{z}_{l}\right)\frac{\partial  }{\partial z_{l} }\left(A_{l}(z,\overline{z})\right)\\&+\left(\overline{z}_{l}\left(1+4\lambda_{l}^{2}\right)+4\lambda_{l}z_{l}\right)\frac{\partial  }{\partial \overline{z}_{l} }\left(A_{l}(z,\overline{z})\right)+ \left(1+4\lambda_{l}^{2}\right)
 A_{l}(z,\overline{z})
 ,\end{split}  \label{90D}
 \end{equation}  
for all $l=1,\dots,N$.

On the other hand, we   combine (\ref{yy}) and  (\ref{oppv}). We obtain 
  \begin{equation}\begin{split}& \quad\quad\quad\quad \quad\quad    \tr_{l} \left(A_{l}(z,\overline{z})\right)q_{l}(z,\overline{z})=   \left(z_{l}\overline{z}_{l} +\lambda_{l}\left(z_{l}^{2}+\overline{z}_{l}^{2}\right) \right)\left( \displaystyle\sum_{I,J\in\mathbb{N}^{N}\atop\left|I\right|+\left|J\right|=p-2} a_{I;J}^{\left(l\right)}i_{l}j_{l} z_{1}^{i_{1}}\dots z_{l}^{i_{l}-1} \dots z_{N}^{i_{N}} \overline{z}_{1}^{j_{1}}\dots \overline{z}_{l}^{j_{l}-1}\dots \overline{z}_{N}^{j_{N}}+\right.\\&\left. \displaystyle\sum_{I,J\in\mathbb{N}^{N}\atop\left|I\right|+\left|J\right|=p-2} a_{I;J}^{\left(l\right)}\lambda_{l} i_{l}\left(i_{l}-1\right) z_{1}^{i_{1}}\dots z_{l}^{i_{l}-2}\dots z_{N}^{i_{N}}\overline{z}_{1}^{j_{1}}\dots  \overline{z}_{N}^{j_{N}}+\displaystyle\sum_{I,J\in\mathbb{N}^{N}\atop\left|I\right|+\left|J\right|=p-2} a_{I;J}^{\left(l\right)}\lambda_{l} j_{l}\left(j_{l}-1\right) z_{1}^{i_{1}}\dots  z_{N}^{i_{N}} \overline{z}_{1}^{j_{1}}\dots \overline{z}_{l}^{j_{l}-2}\dots \overline{z}_{N}^{j_{N}}\right) .\end{split}
\label{B1}\end{equation}
 for all $l=1,\dots,N$.
 
In order to derive a system of equations from  (\ref{90D}), we introduce by (\ref{yy})  the following vectors
\begin{equation} \begin{split}&
X_{1}=X_{1}\left[I;J\right]=\begin{pmatrix}a_{\left(i_{1}+1,\dots,i_{l}-1,\dots,i_{N};j_{1}+1,\dots,j_{l}-1,\dots,j_{N}\right)}^{\left(1\right)} \\ \vdots \\ a_{\left(i_{1},\dots,i_{l},\dots,i_{N};j_{1},\dots,j_{l},\dots,j_{N}\right)}^{\left(l\right)} \\ \vdots \\  a_{\left(i_{1},\dots,i_{l}-1,\dots,i_{N}+1;j_{1},\dots,j_{l}-1,\dots,j_{N}+1\right)}^{\left(N\right)},\end{pmatrix},\\& \quad\quad\quad\quad\hspace{0.1 cm}    X_{2}=X_{2}\left[I;J\right]=\begin{pmatrix}a_{\left(i_{1}+2,\dots,i_{l}-2,\dots,i_{N};j_{1},\dots,j_{N}\right)}^{\left(1\right)}\\ \vdots \\ a_{\left(i_{1},\dots,i_{l},\dots,i_{N};j_{1},\dots,j_{N}\right)}^{\left(l\right)} \\ \vdots \\  a_{\left(i_{1},\dots,i_{l}-2,\dots,i_{N}+2;j_{1},\dots,j_{N}\right)}^{\left(N\right)}\end{pmatrix},\quad  X_{3}=X_{3}\left[I;J\right]=\begin{pmatrix}a_{\left(i_{1},\dots,i_{N};j_{1}+2,\dots,j_{l}-2,\dots,j_{N}\right)}^{\left(1\right)} \\ \vdots \\ a_{\left(i_{1},\dots,i_{N};j_{1},\dots,j_{l}\dots,j_{N}\right)}^{\left(l\right)} \\ \vdots \\  a_{\left(i_{1},\dots,i_{N};j_{1},\dots,j_{l}-2\dots,j_{N}+2\right)}^{\left(N\right)}\label{vectori}\end{pmatrix}.
\end{split} \end{equation}

In order to simplify the computations regarding (\ref{90D}), we introduce   the following   matrix
\begin{equation}\Lambda=\begin{pmatrix}\lambda_{1} & 0     &\dots & 0  \\ 0 & \lambda_{2}   &\dots & 0  \\  \vdots & \vdots &\ddots & \vdots \\ 0  & 0  &   \dots & \lambda_{N} 
\end{pmatrix}, \label{Lambda}
\end{equation}
and then we make the following analysis in (\ref{B1}).

The first very consistent sum of homogeneous terms  is obviously multiplied by $\lambda_{l}z_{l}^{2}$ and $z_{l}\overline{z}_{l}$ and $\lambda_{l}\overline{z}_{l}^{2}$, for all $l\in 1,\dots, N$. This sum    generates obviously  by  (\ref{Lambda})  and (\ref{vectori})   the following terms
\begin{equation}\mathcal{A}X_{1},\quad  \left(\Lambda\mathcal{A}\right)X_{1},\quad \left(\Lambda\mathcal{A} \Lambda\right)X_{1}, 
  \label{idi1}
\end{equation}
where we have used by (\ref{yy})  the following matrix
\begin{equation}  \mathcal{A}=\mathcal{A}\left[I;J\right]= 
\begin{pmatrix} \left(i_{1}+1\right)\left(j_{1}+1\right) & \dots & \left(i_{k-1}+1\right)\left(j_{k-1}+1\right)  & i_{k}j_{k} & \left(i_{k+1}+1\right)\left(j_{k+1}+1\right)  &\dots & \left(i_{N}+1\right)\left(j_{N}+1\right)  \\    \vdots & \ddots & \vdots & \vdots & \vdots & \ddots & \vdots    \\ \left(i_{1}+1\right)\left(j_{1}+1\right) & \dots & \left(i_{k-1}+1\right)\left(j_{k-1}+1\right)  & i_{k}j_{k} & \left(i_{k+1}+1\right)\left(j_{k+1}+1\right)  &\dots & \left(i_{N}+1\right)\left(j_{N}+1\right) \\ \vdots & \ddots & \vdots & \vdots & \vdots & \ddots & \vdots  \\ \left(i_{1}+1\right)\left(j_{1}+1\right) & \dots & \left(i_{k-1}+1\right)\left(j_{k-1}+1\right)  & i_{k}j_{k} & \left(i_{k+1}+1\right)\left(j_{k+1}+1\right)  &\dots & \left(i_{N}+1\right)\left(j_{N}+1\right)
\end{pmatrix}. \label{idiot1}
\end{equation}

The second very consistent sum of homogeneous terms  is obviously multiplied by $\lambda_{l}z_{l}^{2}$ and $z_{l}\overline{z}_{l}$ and $\lambda_{l}\overline{z}_{l}^{2}$, for all $l\in 1,\dots, N$. This sum       generates obviously  by  (\ref{Lambda})  and (\ref{vectori})  the following terms
\begin{equation}
\left( \mathcal{A'}\Lambda\right)X_{2} ,\quad \left(\Lambda\mathcal{A'}\Lambda\right)X_{2} ,\quad \left(\Lambda\mathcal{A'}\Lambda\right)X_{2}, \label{idi2}
\end{equation} 
where we have used by (\ref{yy})  the following matrix
\begin{equation}\mathcal{A'}=\mathcal{A'}\left[I;J\right]= 
\begin{pmatrix} \left(i_{1}+1\right)\left(i_{1}+2\right) & \dots   & \left(i_{k}-1\right)i_{k} & \dots & \left(i_{N}+1\right)\left(i_{N}+2\right) \\    \vdots & \ddots &   \vdots &   \ddots & \vdots    \\ \left(i_{1}+1\right)\left(i_{1}+2\right) & \dots & \left(i_{k}-1\right)i_{k}  & \dots & \left(i_{N}+1\right)\left(i_{N}+2\right) \\ \vdots & \ddots &   \vdots &  \ddots & \vdots  \\ \left(i_{1}+1\right)\left(i_{1}+2\right) & \dots &   \left(i_{k}-1\right)i_{k} &\dots & \left(i_{N}+1\right)\left(i_{N}+2\right)
\end{pmatrix}.\label{idiot2} 
\end{equation}

The third very consistent sum of homogeneous terms  is obviously multiplied by $\lambda_{l}z_{l}^{2}$ and $z_{l}\overline{z}_{l}$ and $\lambda_{l}\overline{z}_{l}^{2}$, for all $l\in 1,\dots, N$. This sum       generates obviously   by   (\ref{Lambda})  and (\ref{vectori})  the following terms
\begin{equation}
\left(\ \mathcal{A''}\Lambda\right)X_{3} ,\quad \left(\Lambda\mathcal{A''}\Lambda\right)X_{3} ,\quad \left(\Lambda\mathcal{A''}\Lambda\right)X_{3}, \label{idi3}
\end{equation}
where we have used by (\ref{yy})  the following matrix
\begin{equation} \mathcal{A''}=\mathcal{A''}\left[I;J\right]= 
\begin{pmatrix} \left(j_{1}+1\right)\left(j_{1}+2\right) & \dots   & \left(j_{k}-1\right)j_{k}  & \dots & \left(j_{N}+1\right)\left(j_{N}+2\right) \\    \vdots & \ddots &   \vdots &   \ddots & \vdots    \\ \left(j_{1}+1\right)\left(j_{1}+2\right) & \dots &  \left(j_{k}-1\right)j_{k}   &\dots & \left(j_{N}+1\right)\left(j_{N}+2\right) \\ \vdots & \ddots &   \vdots &  \ddots & \vdots  \\ \left(j_{1}+1\right)\left(j_{1}+2\right) & \dots &   \left(j_{k}-1\right)j_{k}  &  \dots & \left(j_{N}+1\right)\left(j_{N}+2\right)
\end{pmatrix}. \label{idiot3}
\end{equation}

 It is obvious that   $\Lambda$ does not generally commute with any of the matrices $\mathcal{A}$, $\mathcal{A'}$,  $\mathcal{A''}$ from (\ref{idiot1}), (\ref{idiot2}) and (\ref{idiot3}),  but (\ref{idiot1}), (\ref{idiot2}) and (\ref{idiot3}) provide a better  understanding of the very complicated interactions of of homogeneous terms in (\ref{90D}) as follows.

\subsection{Systems of Equations}  Now, we impose the standard lexicografic order corresponding to \begin{equation}\left(z_{1},z_{2},\dots,z_{N}, \overline{z}_{1},\overline{z}_{2},\dots,\overline{z}_{N}\right)\in\mathbb{C}^{N}\times\mathbb{C}^{N},
\label{ODIN}
\end{equation}   considering     the following vectors   
\begin{equation}\begin{split}& Y_{1}^{t}=\left\{\left(a_{I;0}^{\left(1\right)}\right),\left(a_{I;0}^{\left(2\right)}\right),\dots,\left(a_{I;0}^{\left(N\right)}\right)\right\}_{I\in\mathbb{N}^{N}\atop \left|I\right|=p},\hspace{0.1 cm} Y_{2}^{t}=\left\{ \left(a_{I;0}^{\left(1\right)}\right),\left(a_{I;0}^{\left(2\right)}\right),\dots,\left(a_{I;0}^{\left(N\right)}\right)\right\}_{I,J\in\mathbb{N}^{N}\atop \left|I\right|=p-1, \left|J\right|=1},\dots,\\&  Y_{k}^{t}=\left\{\left(a_{I;0}^{\left(1\right)}\right),\left(a_{I;0}^{\left(2\right)}\right),\dots,\left(a_{I;0}^{\left(N\right)}\right)\right\}_{I,J\in\mathbb{N}^{N}\atop \left|I\right|=p-k+1, \left|J\right|=k-1},\dots, Y_{p+1}^{t}=\left\{\left(a_{I;0}^{\left(1\right)}\right),\left(a_{I;0}^{\left(2\right)}\right),\dots,\left(a_{I;0}^{\left(N\right)}\right)\right\}_{J\in\mathbb{N}^{N}\atop   \left|J\right|=p}.\end{split}
\label{91A}
\end{equation}

 More precisely, we construct very non-trivial  systems of equations   extracting homogeneous terms in (\ref{90D})  using (\ref{op}).   Defining 
 \begin{equation}  V_{1} =\begin{pmatrix}   0\\  \vdots \\0 \\   \lambda_{1}\tilde{i}_{1}\left(\tilde{i}_{1}-1\right)  \\0\\ \dots \\0 \\ \lambda_{l}\tilde{i}_{l}\left(\tilde{i}_{l}-1\right) \\0\\ \vdots \\ \lambda_{N}\tilde{i}_{N}\left(\tilde{i}_{N}-1\right) \\0\\ \vdots  \\ 0 
\end{pmatrix}
,\quad V_{0} =\begin{pmatrix}  0\\ \vdots \\0\\0\\0\\  \vdots\\  0\\0\\0\\  \vdots \\0\\0\\0\\  \vdots  
\\ 0\end{pmatrix}, \label{special1}
\end{equation}
we obtain    the following  system of equations
 \begin{equation}\begin{pmatrix} \mathcal{M}_{1,1} & \mathcal{M}_{1,2} & \mathcal{M}_{1,3} & \mbox{O}_{N^{p}} & \mbox{O}_{N^{p}} &\dots & \mbox{O}_{N^{p}} & \mbox{O}_{N^{p}} & \mbox{O}_{N^{p}} & \mbox{O}_{N^{p}} \\ \mathcal{M}_{2,1} & \mathcal{M}_{2,2} & \mathcal{M}_{2,3} & \mathcal{M}_{2,4} & \mbox{O}_{N^{p}} &\dots & \mbox{O}_{N^{p}} & \mbox{O}_{N^{p}} & \mbox{O}_{N^{p}} & \mbox{O}_{N^{p}} \\ \mathcal{M}_{3,1}  & \mathcal{M}_{3,2}  & \mathcal{M}_{3,3}  & \mathcal{M}_{3,4}  & \mathcal{M}_{3,5}  &\dots & \mbox{O}_{N^{p}} & \mbox{O}_{N^{p}} & \mbox{O}_{N^{p}} & \mbox{O}_{N^{p}} \\ \mbox{O}_{N^{p}} & \mathcal{M}_{4,2}  & \mathcal{M}_{4,3}  & \mathcal{M}_{4,4}  & \mathcal{M}_{4,5}  &\dots & \mbox{O}_{N^{p}} & \mbox{O}_{N^{p}} & \mbox{O}_{N^{p}}& \mbox{O}_{N^{p}} \\ \mbox{O}_{N^{p}} & \mbox{O}_{N^{p}} & \mathcal{M}_{5,3} & \mathcal{M}_{5,4}& \mathcal{M}_{5,5} &\dots & \mbox{O}_{N^{p}} &\mbox{O}_{N^{p}} & \mbox{O}_{N^{p}} & \mbox{O}_{N^{p}} \cr \mbox{O}_{N^{p}} & \mbox{O}_{N^{p}} & \mbox{O}_{N^{p}} & \mathcal{M}_{6,4} & \mathcal{M}_{6,5} &\dots & \mbox{O}_{N^{p}} & \mbox{O}_{N^{p}} & \mbox{O}_{N^{p}} & \mbox{O}_{N^{p}}  \\ \vdots & \vdots & \vdots & \vdots & \vdots & \ddots & \vdots & \vdots & \vdots & \vdots \\ \mbox{O}_{N^{p}} & \mbox{O}_{N^{p}} & \mbox{O}_{N^{p}} & \mbox{O}_{N^{p}} & \mbox{O}_{N^{p}} &\dots & \mathcal{M}_{p-2,p-2}  & \mathcal{M}_{p-2,p-1}  & \mathcal{M}_{p-2,p} & \mbox{O}_{N^{p}} \\ \mbox{O}_{N^{p}} & \mbox{O}_{N^{p}} & \mbox{O}_{N^{p}} & \mbox{O}_{N^{p}} & \mbox{O}_{N^{p}} &\dots & \mathcal{M}_{p-1,p-2}  & \mathcal{M}_{p-1,p-1}  &\mathcal{M}_{p-1,p}  & \mathcal{M}_{p-1,p+1}  \\ \mbox{O}_{N^{p}} & \mbox{O}_{N^{p}} & \mbox{O}_{N^{p}} & \mbox{O}_{N^{p}} & \mbox{O}_{N^{p}} &\dots & \mathcal{M}_{p,p-2}  & \mathcal{M}_{p,p-1}  & \mathcal{M}_{p,p}  & \mathcal{M}_{p,p+1}  \\ \mbox{O}_{N^{p}} & \mbox{O}_{N^{p}} & \mbox{O}_{N^{p}} & \mbox{O}_{N^{p}} &\mbox{O}_{N^{p}} &\dots & \mbox{O}_{N^{p}} & \mathcal{M}_{p+1,p-1}  & \mathcal{M}_{p+1,p}  & \mathcal{M}_{p+1,p+1}           
\end{pmatrix}\begin{pmatrix} Y_{1}\\Y_{2}\\Y_{3}\\Y_{4}\\Y_{5}\\ \vdots\\ Y_{p-3}\\ Y_{p-2} \\ Y_{p-1}\\ Y_{p}\\ Y_{p+1}
\end{pmatrix}=\begin{pmatrix} V_{1} \\ V_{0} \\ V_{0} \\ V_{0}\\ V_{0}  \\   \vdots \\ V_{0} \\ V_{0}\\ V_{0}\\ V_{0}\\ V_{0}
\end{pmatrix}, \label{beb1}
\end{equation}
where its   elements are matrices, which  are defined as follows
\begin{equation} \begin{split}& \mathcal{M}_{k-2,k}=\mathcal{W}_{k}^{'},\quad\quad\quad\quad  \quad\quad \quad\quad\quad \hspace{0.1 cm}  \mbox{for all $k=3,\dots,p+1$,} \\&\mathcal{M}_{k-1,k}=\mathcal{V}_{k}^{(0)}+\mathcal{W}_{k}+\mathcal{O}_{k}^{'} ,  \quad\quad\quad\hspace{0.2 cm}\mbox{for all $k=2,\dots,p+1 $,}\\&\quad \mathcal{M}_{k,k}= \mathcal{O}_{k}^{(0)}+\mathcal{O}_{k}+\mathcal{V}_{k}^{'}+\mathcal{W}_{k}^{''}  ,\quad\mbox{for all $k=1,\dots,p+1$,}  \\& \mathcal{M}_{k+1,k}=\mathcal{W}_{k}^{(0)}+\mathcal{V}_{k}+\mathcal{O}_{k}^{''} ,  \quad\quad\hspace{0.19  cm}\quad\mbox{for all $k=1,\dots,p$,} \\&  \mathcal{M}_{k+2,k}=\mathcal{V}_{k}^{''},\quad\quad\quad\quad\quad\quad\quad\quad  \quad \hspace{0.21 cm}  \mbox{for all $k=1,\dots,p-1$.} \end{split}  \label{lil}
\end{equation} 

These occurring  matrices   are defined as follows:

$\bf{The\hspace{0.1 cm}matrices:}$
$$\bf{\left\{\mathcal{O}_{k}^{(0)}\right\}_{k=1,\dots,p+1}\hspace{0.1 cm}, \left\{\mathcal{V}_{k}^{(0)}\right\}_{k=2,\dots,p+1},\hspace{0.1 cm} \left\{\mathcal{W}_{k}^{(0)}\right\}_{k=1,\dots,p}}\in \mathcal{M}_{N^{p}\times N^{p}}\left(\mathbb{C}\right).$$

 Occurring from (\ref{90D}), we have used in (\ref{lil}) the following matrices 
\begin{equation}\mathcal{O}_{k}^{(0)}\left[I;J\right]=\begin{pmatrix} 1& \dots & 0 & \dots & 0 & \dots & 0\\ \vdots & \ddots& \vdots & \ddots & \vdots& \ddots & \vdots \\ 0&\dots & \left(i_{1}+ j_{1}+1\right)\left(1+4\lambda_{1}^{2}\right)    &\dots & 0   &   \dots & 0 &  \\ \vdots & \ddots& \vdots & \ddots & \vdots& \ddots & \vdots  \\ 0&\dots & 0& \dots &  \left(i_{N}+ j_{N}+1\right)\left(1+4\lambda_{N}^{2}\right)    & \dots  & 0  \\   \vdots & \ddots& \vdots & \ddots & \vdots& \ddots & \vdots\\ 0&\dots &0 &  \dots & 0 & \dots &  1  \end{pmatrix}, \label{extra1}
\end{equation} 
such that   
\begin{equation}j_{1}+j_{2}+\dots+j_{N}+i_{1}+i_{2}+\dots+i_{N}=p ,\quad j_{1}+j_{2}+\dots+j_{N}=k-1.\label{lili}
\end{equation}

Then (\ref{lili}) defines in (\ref{lil}) also the following   matrices
\begin{equation} \mathcal{V}_{k}^{(0)}\left[I;J\right]=\begin{pmatrix} 1& \dots & 0 & \dots & 0 & \dots & 0\\ \vdots & \ddots& \vdots & \ddots & \vdots& \ddots & \vdots \\ 0&\dots & 4\lambda_{1} i_{1}    &\dots & 0   &   \dots & 0 &  \\ \vdots & \ddots& \vdots & \ddots & \vdots& \ddots & \vdots  \\ 0&\dots & 0& \dots & 4\lambda_{N} i_{N}     & \dots  & 0  \\   \vdots & \ddots& \vdots & \ddots & \vdots& \ddots & \vdots\\ 0&\dots &0 &  \dots & 0 & \dots &  1  \end{pmatrix},\quad\mbox{for all $k=2,\dots,p+1 $,}\label{extra2} 
\end{equation} 
and respectively, the following matrices
\begin{equation} \mathcal{W}_{k}^{(0)}\left[I;J\right]=\begin{pmatrix} 1& \dots & 0 & \dots & 0 & \dots & 0\\ \vdots & \ddots& \vdots & \ddots & \vdots& \ddots & \vdots \\ 0&\dots & 4\lambda_{1} j_{1}    &\dots & 0   &   \dots & 0 &  \\ \vdots & \ddots& \vdots & \ddots & \vdots& \ddots & \vdots  \\ 0&\dots & 0& \dots & 4\lambda_{N} j_{N}     & \dots  & 0  \\   \vdots & \ddots& \vdots & \ddots & \vdots& \ddots & \vdots\\ 0&\dots &0 &  \dots & 0 & \dots &  1  \end{pmatrix},\quad\mbox{for all $k=1,\dots,p $,}\label{extra3} 
\end{equation}

Thus  in (\ref{lil}), the following matrices are by (\ref{extra1}), (\ref{extra2}) and (\ref{extra3})  defined as follows
\begin{equation} \begin{split}&\hspace{0.05 cm}
\mathcal{O}_{k}^{(0)}=\mathcal{O}_{k}^{(0)}\left[p-2,\dots,0;0,\dots,0\right]\dots\mathcal{O}_{k}^{(0)}\left[p-3,\dots,0;1,\dots,0\right]\dots\mathcal{O}_{k}^{(0)}\left[0,\dots,0;0,\dots,p-2\right],\quad\hspace{0.22 cm}\mbox{for all $k=1,\dots,p+1$,}  \\&\hspace{0.1 cm}\mathcal{V}_{k}^{(0)}=\mathcal{V}_{k}^{(0)}\left[p-2,\dots,0;0,\dots,0\right]\dots\mathcal{V}_{k}^{(0)}\left[p-3,\dots,0;1,\dots,0\right]\dots \mathcal{V}_{k}^{(0)}\left[0,\dots,0;0,\dots,p-2\right],\quad\quad \mbox{for all $k=2,\dots,p+1$,}  \\& \mathcal{W}_{k}^{(0)}=\mathcal{W}_{k}^{(0)}\left[p-2,\dots,0;0,\dots,0\right]\dots\mathcal{W}_{k}^{(0)}\left[p-3,\dots,0;1,\dots,0\right]\dots \mathcal{W}_{k}^{(0)}\left[0,\dots,0;0,\dots,p-2\right],\quad\mbox{for all $k=1,\dots,p$.} \end{split} \label{extra4}
\end{equation}

$\bf{The\hspace{0.1 cm}matrices:}$
$$\bf{\left\{\mathcal{O}_{k}\right\}_{k=1,\dots,p+1},\hspace{0.1 cm}  \left\{\mathcal{V}_{k}\right\}_{k=1,\dots,p},\hspace{0.1 cm} \left\{\mathcal{W}_{k}\right\}_{k=2,\dots,p+1}}\in \mathcal{M}_{N^{p}\times N^{p}}\left(\mathbb{C}\right).$$

 The first   matrix from (\ref{idiot1}) induces by (\ref{yy}), (\ref{B1}) and   (\ref{lili})  the following matrices  
\begin{equation}   \mathcal{O}_{k}\left[I;J\right]=
\begin{pmatrix}  1 & \dots   & 0 & \dots & 0 & \dots & 0& \dots & 0  \\ \vdots & \ddots & \vdots & \ddots & \vdots & \ddots & \vdots& \ddots & \vdots   \\ 0 & \dots & \left(i_{1}+1\right)\left(j_{1}+1\right) & \dots & i_{k}j_{k} & \dots & \left(i_{N}+1\right)\left(j_{N}+1\right) & \dots & 0 
\\ \vdots & \ddots & \vdots  & \ddots & \vdots  & \ddots & \vdots  & \ddots & \vdots  \\ 0 & \dots & \left(i_{1}+1\right)\left(j_{1}+1\right) & \dots & i_{k}j_{k} & \dots & \left(i_{N}+1\right)\left(j_{N}+1\right) & \dots & 0 \\ \vdots  & \ddots & \vdots  & \ddots & \vdots  & \ddots & \vdots  & \ddots & \vdots  \\ 0 & \dots & \left(i_{1}+1\right)\left(j_{1}+1\right) & \dots & i_{k}j_{k} & \dots & \left(i_{N}+1\right)\left(j_{N}+1\right) & \dots & 0 \\ \vdots  & \ddots & \vdots  & \ddots & \vdots  & \ddots & \vdots  & \ddots & \vdots  \\ 0  & \dots & 0  & \dots &0  & \dots & 0  & \dots & 1  
\end{pmatrix},  \label{extra1A}
\end{equation}
for all $k=1,\dots,p+1$.  

This matrix (\ref{extra1A}) has the characteristic that $i_{k}j_{k}$    is  the diagonal entry on the following row $$\left(i_{1},\dots,i_{N};j_{1},\dots,j_{N}\right)\in\mathbb{N}^{N}\times\mathbb{N}^{N},$$  according to the    lexicografic order  related to  (\ref{ODIN}), otherwise having   $1$ as diagonal entry.

Similarly, we consider by (\ref{yy})  and (\ref{idiot1}) the following matrices
 \begin{equation}   \mathcal{V}_{k}\left[I;J\right]=
\begin{pmatrix}  1 & \dots   & 0 & \dots & 0 & \dots & 0& \dots & 0  \\ \vdots & \ddots & \vdots & \ddots & \vdots & \ddots & \vdots& \ddots & \vdots   \\ 0 & \dots & \left(i_{1}+1\right)\left(j_{1}+1\right)
 \lambda_{1} & \dots & i_{k}j_{k}
 \lambda_{k} & \dots & \left(i_{N}+1\right)\left(j_{N}+1\right)
 \lambda_{N} & \dots & 0 
\\ \vdots & \ddots & \vdots  & \ddots & \vdots  & \ddots & \vdots  & \ddots & \vdots  \\ 0 & \dots & \left(i_{1}+1\right)\left(j_{1}+1\right)
 \lambda_{1} & \dots & i_{k}j_{k}
 \lambda_{k} & \dots & \left(i_{N}+1\right)\left(j_{N}+1\right)
 \lambda_{N} & \dots & 0 \\ \vdots  & \ddots & \vdots  & \ddots & \vdots  & \ddots & \vdots  & \ddots & \vdots  \\ 0 & \dots & \left(i_{1}+1\right)\left(j_{1}+1\right)
 \lambda_{1} & \dots & i_{k}j_{k}
 \lambda_{k} & \dots & \left(i_{N}+1\right)\left(j_{N}+1\right)
 \lambda_{N} & \dots & 0 \\ \vdots  & \ddots & \vdots  & \ddots & \vdots  & \ddots & \vdots  & \ddots & \vdots  \\ 0  & \dots & 0  & \dots &0  & \dots & 0  & \dots & 1  
\end{pmatrix},  \label{extra2A}
\end{equation}
for all $k=1,\dots,p$.  

This matrix  (\ref{extra2A}) has the characteristic that $i_{k}j_{k}\lambda_{k}$    is  the diagonal entry  on the  following row $$ \left(i_{1},\dots,i_{k}+1,\dots,i_{N};j_{1},\dots,j_{k}-1,\dots,j_{N}\right)\in\mathbb{N}^{N}\times\mathbb{N}^{N},$$  according to the corresponding  lexicografic order  related to  (\ref{ODIN}), otherwise having   $1$ as diagonal entry. 

Then, (\ref{extra2A}) induces by (\ref{yy})  and (\ref{idiot1}) analogously another matrices denoted  as follows
\begin{equation}\mathcal{W}_{k}\left[I;J\right],\quad\mbox{for all $k=2,\dots,p+1$,}\label{extra3A}\end{equation} having the characteristic that $i_{k}j_{k}\lambda_{k}$    is  the diagonal entry  on the  following row   $$ \left(i_{1},\dots,i_{k}-1,\dots,i_{N};j_{1},\dots,j_{k}+1,\dots,j_{N}\right)\in\mathbb{N}^{N}\times\mathbb{N}^{N},$$ according to the    lexicografic order  related to  (\ref{ODIN}), otherwise having   $1$ as diagonal entry.  

 Thus in (\ref{lil}),  the following matrices are by (\ref{extra1A}), (\ref{extra2A}) and (\ref{extra3A})  naturally defined  as follows
\begin{equation} \begin{split}&\hspace{0.05 cm}
\mathcal{O}_{k}=\mathcal{O}_{k}\left[p-2,\dots,0;0,\dots,0\right]\dots\mathcal{O}_{k}\left[p-3,\dots,0;1,\dots,0\right]\dots \mathcal{O}_{k}\left[0,\dots,0;0,\dots,p-2\right], \hspace{0.16 cm}\quad\mbox{for all $k=1,\dots,p+1$,} \\& \hspace{0.1 cm}\mathcal{V}_{k}=\mathcal{V}_{k}\left[p-2,\dots,0;0,\dots,0\right]\dots\mathcal{V}_{k}\left[p-3,\dots,0;1,\dots,0\right]\dots \mathcal{V}_{k}\left[0,\dots,0;0,\dots,p-2\right], \quad\quad\mbox{for all $k=1,\dots,p$,}  \\& \mathcal{W}_{k}=\mathcal{W}_{k}\left[p-2,\dots,0;0,\dots,0\right]\dots\mathcal{W}_{k}\left[p-3,\dots,0;1,\dots,0\right]\dots \mathcal{W}_{k}\left[0,\dots,0;0,\dots,p-2\right],\quad\mbox{for all $k=2,\dots,p+1$.}\end{split} \label{extra4A}
\end{equation}

It is important to better explain these products from (\ref{extra4}) and (\ref{extra4A}). These matrices are commuting, because there are no commune non-trivial rows. Thus, (\ref{extra4}) and (\ref{extra4A}) are clear because   (\ref{extra1}),(\ref{extra2}),(\ref{extra3}), (\ref{extra1A}), (\ref{extra2A}) and (\ref{extra3A}) define classes of commuting matrices.

$\bf{The\hspace{0.1 cm}matrices:}$
$$\bf{\left\{\mathcal{O}_{k}^{'}\right\}_{k=2,\dots,p+1},\hspace{0.1 cm}  \left\{\mathcal{V}_{k}^{'}\right\}_{k=1,\dots,p+1},\hspace{0.1 cm} \left\{\mathcal{W}_{k}^{'}\right\}_{k=3,\dots,p+1}}\in \mathcal{M}_{N^{p}\times N^{p}}\left(\mathbb{C}\right).$$ 

The first  matrix from (\ref{idiot2}) induces   by   (\ref{yy}),   (\ref{B1}) and (\ref{lili})    the following matrices 
 \begin{equation}      \mathcal{O'}_{k}\left[I;J\right]= \begin{pmatrix}  1 & \dots   & 0 & \dots & 0 & \dots & 0& \dots & 0  \\ \vdots & \ddots & \vdots & \ddots & \vdots & \ddots & \vdots& \ddots & \vdots   \\ 0 & \dots & \left(i_{1}+1\right)\left(i_{1}+2\right)\lambda_{1}  & \dots & \left(i_{k}-1\right)i_{k}\lambda_{k} & \dots & \left(i_{N}+1\right)\left(i_{N}+2\right)\lambda_{N} & \dots & 0 
\\ \vdots & \ddots & \vdots  & \ddots & \vdots  & \ddots & \vdots  & \ddots & \vdots  \\ 0 & \dots & \left(i_{1}+1\right)\left(i_{1}+2\right)\lambda_{1} & \dots & \left(i_{k}-1\right)i_{k}\lambda_{k} & \dots & \left(i_{N}+1\right)\left(i_{N}+2\right)\lambda_{N} & \dots & 0 \\ \vdots  & \ddots & \vdots  & \ddots & \vdots  & \ddots & \vdots  & \ddots & \vdots  \\ 0 & \dots & \left(i_{1}+1\right)\left(i_{1}+2\right)\lambda_{1}  & \dots & \left(i_{k}-1\right)i_{k}\lambda_{k} & \dots & \left(i_{N}+1\right)\left(i_{N}+2\right)\lambda_{N} & \dots & 0 \\ \vdots  & \ddots & \vdots  & \ddots & \vdots  & \ddots & \vdots  & \ddots & \vdots  \\ 0  & \dots & 0  & \dots &0  & \dots & 0  & \dots & 1  
\end{pmatrix},\label{extra1B} 
\end{equation} 
for all  $k=2,\dots,p+1$.  

This matrix   (\ref{extra1B}) has the characteristic that  $\left(i_{k}-1\right)i_{k}\lambda_{k}$ is  the diagonal entry  of the following  row  $$\left(i_{1},\dots,i_{k}-1,\dots,i_{N};j_{1},\dots,j_{k}+1,\dots,j_{N}\right)\in\mathbb{N}^{N}\times\mathbb{N}^{N},$$  according to the  lexicografic order  related to  (\ref{ODIN}), otherwise having   $1$ as diagonal entry. 

Then, the second matrix from (\ref{idiot2}) induces   by (\ref{yy})   the following matrices
  \begin{equation}     \mathcal{V'}_{k}\left[I;J\right]= \begin{pmatrix}  1 & \dots   & 0 & \dots & 0 & \dots & 0& \dots & 0  \\ \vdots & \ddots & \vdots & \ddots & \vdots & \ddots & \vdots& \ddots & \vdots   \\ 0 & \dots & \left(i_{1}+1\right)\left(i_{1}+2\right)\lambda_{1}^{2} & \dots & \left(i_{k}-1\right)i_{k}\lambda_{1} \lambda_{k}   & \dots & \left(i_{N}+1\right)\left(i_{N}+2\right)\lambda_{1}\lambda_{N}   & \dots & 0 
\\ \vdots & \ddots & \vdots  & \ddots & \vdots  & \ddots & \vdots  & \ddots & \vdots  \\ 0 & \dots & \left(i_{1}+1\right)\left(i_{1}+2\right)\lambda_{1}\lambda_{k}   & \dots & \left(i_{k}-1\right)i_{k}\lambda_{k}^{2} & \dots & \left(i_{N}+1\right)\left(i_{N}+2\right)\lambda_{N}\lambda_{k}  & \dots & 0 \\ \vdots  & \ddots & \vdots  & \ddots & \vdots  & \ddots & \vdots  & \ddots & \vdots  \\ 0 & \dots & \left(i_{1}+1\right)\left(i_{1}+2\right)\lambda_{1}\lambda_{N}   & \dots & \left(i_{k}-1\right)i_{k}\lambda_{k}\lambda_{N}  & \dots & \left(i_{N}+1\right)\left(i_{N}+2\right)\lambda_{N}^{2} & \dots & 0 \\ \vdots  & \ddots & \vdots  & \ddots & \vdots  & \ddots & \vdots  & \ddots & \vdots  \\ 0  & \dots & 0  & \dots &0  & \dots & 0  & \dots & 1  
\end{pmatrix},\label{extra2B} 
\end{equation} 
for all  $k=1,\dots,p+1$.  

This matrix (\ref{extra2B}) has the characteristic that  $\left(i_{k}-1\right)i_{k}\lambda_{k}^{2}$ is  the diagonal entry  on the  following row   $$\left(i_{1}, \dots,i_{N};j_{1},\dots,j_{N}\right)\in\mathbb{N}^{N}\times\mathbb{N}^{N},$$  according to the   lexicografic order  related to  (\ref{ODIN}), otherwise having   $1$ as diagonal entry.  

Also the matrix (\ref{extra2B}) induces by (\ref{yy})  and (\ref{idiot2}) another matrices denoted   as follows 
\begin{equation}\mathcal{W'}_{k}\left[I;J\right],\quad\mbox{for all $k=3,\dots,p+1$.}\label{extra3B}\end{equation}

 This matrix (\ref{extra3B}) has the characteristic that $\left(i_{k}-1\right) i_{k}  \lambda_{k}^{2}$ is  the diagonal entry  on the  following row  $$\left(i_{1},\dots,i_{k}-2,\dots,i_{N};j_{1},\dots,j_{k}+2, \dots,j_{N}\right)\in\mathbb{N}^{N}\times\mathbb{N}^{N},$$  according to the    lexicografic order  related to  (\ref{ODIN}), otherwise having   $1$ as diagonal entry. 
 
Thus in (\ref{lil}), the following matrices are by (\ref{extra1B}), (\ref{extra2B}) and (\ref{extra3B})  naturally defined as follows
\begin{equation} \begin{split}&\hspace{0.05 cm}
\mathcal{O'}_{k}=\mathcal{O'}_{k}\left[p-2,\dots,0;0,\dots,0\right]\dots\mathcal{O'}_{k}\left[p-3,\dots,0;1,\dots,0\right]\dots\mathcal{O'}_{k}\left[0,\dots,0;0,\dots,p-2\right], \quad\hspace{0.19 cm}\mbox{for all $k=2,\dots,p+1$,} \\& \hspace{0.1 cm}\mathcal{V'}_{k}=\mathcal{V'}_{k}\left[p-2,\dots,0;0,\dots,0\right]\dots\mathcal{V'}_{k}\left[p-3,\dots,0;1,\dots,0\right]\dots\mathcal{V'}_{k}\left[0,\dots,0;0,\dots,p-2\right],\quad\quad\mbox{for all $k=1,\dots,p+1$,}  \\& \mathcal{W'}_{k}=\mathcal{W'}_{k}\left[p-2,\dots,0;0,\dots,0\right]\dots\mathcal{W'}_{k}\left[p-3,\dots,0;1,\dots,0\right]\dots\mathcal{W'}_{k}\left[0,\dots,0;0,\dots,p-2\right],\quad\mbox{for all $k=3,\dots,p+1$.}\end{split} \label{extra4B}
\end{equation}

$\bf{The\hspace{0.1 cm}matrices:}$
$$ \bf{\left\{\mathcal{O}_{k}^{''}\right\}_{k=1,\dots,p},\hspace{0.1 cm}  \left\{\mathcal{V}_{k}^{''}\right\}_{k=1,\dots,p-1},\hspace{0.1 cm} \left\{\mathcal{W}_{k}^{''}\right\}_{k=1,\dots,p+1}}\in \mathcal{M}_{N^{p}\times N^{p}}\left(\mathbb{C}\right). $$

The first  matrix from (\ref{idiot3}) induces by    (\ref{yy}), (\ref{B1}) and  (\ref{lili}) the following matrices 
 \begin{equation}  \mathcal{O''}_{k}\left[I;J\right]= \begin{pmatrix}  1 & \dots   & 0 & \dots & 0 & \dots & 0& \dots & 0  \\ \vdots & \ddots & \vdots & \ddots & \vdots & \ddots & \vdots& \ddots & \vdots   \\ 0 & \dots & \left(j_{1}+1\right)\left(j_{1}+2\right)\lambda_{1}  & \dots &  \left(j_{k}-1\right)j_{k}\lambda_{k} & \dots & \left(j_{N}+1\right)\left(j_{N}+2\right)\lambda_{N} & \dots & 0 
\\ \vdots & \ddots & \vdots  & \ddots & \vdots  & \ddots & \vdots  & \ddots & \vdots  \\ 0 & \dots & \left(j_{1}+1\right)\left(j_{1}+2\right)\lambda_{1} & \dots &  \left(j_{k}-1\right)j_{k}\lambda_{k} & \dots & \left(j_{N}+1\right)\left(j_{N}+2\right)\lambda_{N} & \dots & 0 \\ \vdots  & \ddots & \vdots  & \ddots & \vdots  & \ddots & \vdots  & \ddots & \vdots  \\ 0 & \dots & \left(j_{1}+1\right)\left(j_{1}+2\right)\lambda_{1}  & \dots &  \left(j_{k}-1\right)j_{k}\lambda_{k} & \dots & \left(j_{N}+1\right)\left(j_{N}+2\right)\lambda_{N} & \dots & 0 \\ \vdots  & \ddots & \vdots  & \ddots & \vdots  & \ddots & \vdots  & \ddots & \vdots  \\ 0  & \dots & 0  & \dots &0  & \dots & 0  & \dots & 1  
\end{pmatrix} , \label{extra1C}
\end{equation}
for all $k=1,\dots,p$.  
 
This matrix   (\ref{extra1C}) has the characteristic that $\left(j_{k}-1\right)j_{k}\lambda_{k}$    is  the diagonal entry  on the  following row $$\left(i_{1},\dots,i_{k}+1,\dots,i_{N};j_{1},\dots,j_{k}-1, \dots,j_{N}\right)\in\mathbb{N}^{N}\times\mathbb{N}^{N},$$    according to the  lexicografic order  related to  (\ref{ODIN}), otherwise having   $1$ as diagonal entry. 

Similarly as previously, we consider by (\ref{yy})   the following matrices
 \begin{equation} \mathcal{V''}_{k}\left[I;J\right]= \begin{pmatrix}  1 & \dots   & 0 & \dots & 0 & \dots & 0& \dots & 0  \\ \vdots & \ddots & \vdots & \ddots & \vdots & \ddots & \vdots& \ddots & \vdots   \\ 0 & \dots & \left(j_{1}+1\right)\left(j_{1}+2\right)\lambda_{1}^{2} & \dots & \left(j_{k}-1\right)j_{k}\lambda_{1}\lambda_{k}  & \dots & \left(j_{N}+1\right)\left(j_{N}+2\right)\lambda_{1}\lambda_{N}  & \dots & 0 
\\ \vdots & \ddots & \vdots  & \ddots & \vdots  & \ddots & \vdots  & \ddots & \vdots  \\ 0 & \dots & \left(j_{1}+1\right)\left(j_{1}+2\right)\lambda_{1}\lambda_{k}   & \dots & \left(j_{k}-1\right)j_{k}\lambda_{k}^{2} & \dots & \left(j_{N}+1\right)\left(j_{N}+2\right)\lambda_{k} \lambda_{N} & \dots & 0 \\ \vdots  & \ddots & \vdots  & \ddots & \vdots  & \ddots & \vdots  & \ddots & \vdots  \\ 0 & \dots & \left(j_{1}+1\right)\left(j_{1}+2\right)\lambda_{1}\lambda_{N}   & \dots & \left(j_{k}-1\right)j_{k}\lambda_{k}\lambda_{N}   & \dots & \left(j_{N}+1\right)\left(j_{N}+2\right)\lambda_{N}^{2} & \dots & 0 \\ \vdots  & \ddots & \vdots  & \ddots & \vdots  & \ddots & \vdots  & \ddots & \vdots  \\ 0  & \dots & 0  & \dots &0  & \dots & 0  & \dots & 1  
\end{pmatrix}  , \label{extra2C}
\end{equation}
for all $k=1,\dots,p-1$.  

This matrix   (\ref{extra2C}) has the characteristic that $\left(j_{k}-1\right)j_{k}\lambda_{k}^{2}$    is  the diagonal entry  on the  following row   
$$\left(i_{1},\dots,i_{k}+2,\dots,i_{N};j_{1},\dots,j_{k}-2, \dots,j_{N}\right)\in\mathbb{N}^{N}\times\mathbb{N}^{N},$$  
 according to the   lexicografic order related to (\ref{ODIN}), otherwise having   $1$ as diagonal entry.
 
Also the matrix  (\ref{extra2C}) induces by (\ref{yy})   another matrices denoted   as follows 
\begin{equation}\mathcal{W''}_{k}\left[I;J\right]\in\mathcal{M}_{N^{p}\times N^{p}}\left(\mathbb{C}\right) ,\quad\mbox{ for all $k=1,\dots,p+1$}.\label{extra3C}\end{equation}  

This matrix (\ref{extra3C}) has the characteristic that $\left(j_{k}-1\right)j_{k} \lambda_{k}^{2}$   is  the diagonal entry  on the  following row 
$$\left(i_{1}, i_{2},\dots,i_{N};j_{1}, j_{2}, \dots,j_{N}\right)\in\mathbb{N}^{N}\times\mathbb{N}^{N},$$ 
according to the   lexicografic order  related to  (\ref{ODIN}), otherwise having  $1$ as diagonal entry.

Thus, in (\ref{lil}) the following matrices are by (\ref{extra1C}), (\ref{extra2C}) and (\ref{extra3C})  defined as follows
\begin{equation} \begin{split}&\hspace{0.05 cm}
\mathcal{O''}_{k}=\mathcal{O''}_{k}\left[p-2,\dots,0;0,\dots,0\right]\dots\mathcal{O''}_{k}\left[p-3,\dots,0;1,\dots,0\right]\dots\mathcal{O''}_{k}\left[0,\dots,0;0,\dots,p-2\right],\quad\hspace{0.1 cm} \hspace{0.1 cm}\mbox{for all $k=1,\dots,p$,}  \\&\hspace{0.1 cm} \mathcal{V''}_{k}=\mathcal{V''}_{k}\left[p-2,\dots,0;0,\dots,0\right]\dots\mathcal{V''}_{k}\left[p-3,\dots,0;1,\dots,0\right]\dots \mathcal{V''}_{k}\left[0,\dots,0;0,\dots,p-2\right],\quad\quad \mbox{for all $k=1,\dots,p-1$,}  \\& \mathcal{W''}_{k}=\mathcal{W''}_{k}\left[p-2,\dots,0;0,\dots,0\right]\dots\mathcal{W''}_{k}\left[p-3,\dots,0;1,\dots,0\right]\dots \mathcal{W''}_{k}\left[0,\dots,0;0,\dots,p-2\right],\quad\mbox{for all $k=1,\dots,p+1$.}\end{split} \label{extra4C}
\end{equation}
 
We were using the following obvious observation 
\begin{equation}   \#\left\{\left(i_{1},i_{2},\dots,i_{N}:j_{1},j_{2},\dots,j_{N}\right)\in \mathbb{N}^{N}\times \mathbb{N}^{N};\quad  \mbox{   $i_{1}+i_{2}+\dots+i_{N}=p-k$, $j_{1}+j_{2}+\dots+j_{N}=k$}\right\}=N^{p},\hspace{0.1 cm} \mbox{for all $k=0,\dots,p$.}\label{set} \end{equation}

The computations continue precisely as in \cite{val8}. We arrive at
\subsection{Fischer Normalization $G$-Spaces\cite{bu4},\cite{bu5}}The Fischer Decomposition from (\ref{op}) gives
\begin{equation}z^{I}=A_{I}(z,\overline{z}) q_{l}(z,\overline{z})+C_{I}(z,\overline{z}),\quad\mbox{where}\hspace{0.1 cm}  \tr_{l}\left(C_{I}(z,\overline{z})\right)=0,\label{opt}
\end{equation}
where $I \in \mathbb{N}$ having   length $p$ and   $l=1,\dots,N$.
 
This set of homogeneous polynomials, derived from (\ref{opt}),   dictates   normalizations based on  certain   Normalization Spaces, which is defined using the generalized version of  the Fischer Decomposition\cite{sh}: for a given  real homogeneous polynomial of degree $p\geq 1$ in $(z,\overline{z})$ denoted by $P(z,\overline{z})$, we  have
\begin{equation} \begin{split}& P(z,\overline{z})=P_{1}(z,\overline{z})q_{l}(z,\overline{z})+R_{1}(z,\overline{z}),\quad\mbox{where $\tr_{l}\left(R_{1}(z,\overline{z})\right)=0$ and:}\\& \quad\quad\quad\quad\hspace{0.1 cm}  R_{1}(z,\overline{z})=\displaystyle\sum_{I\in\mathbb{N}^{N}\atop {\left|I\right|=p}} \left(a_{I}C_{I}(z,\overline{z})+b_{I}\overline{C_{I}(z,\overline{z})}\right)+R_{1,0}(z,\overline{z}),\quad\mbox{where $R_{1,0}(z,\overline{z})\in \displaystyle\bigcap_{I\in\mathbb{N}^{N},\hspace{0.1 cm}   \atop {\left|I\right|=p}} \ker  C^{\star}_{I}  \cap  \ker  \overline{C}^{\star}_{I} $,}
\\& P_{1}(z,\overline{z})=P_{2}(z,\overline{z})q_{l}(z,\overline{z})+R_{2}(z,\overline{z}),\quad\mbox{where $\tr_{l}\left(R_{2}(z,\overline{z})\right)=0$ and:}\\& \quad\quad\quad\quad\hspace{0.1 cm}  R_{2}(z,\overline{z})=\displaystyle\sum_{I\in\mathbb{N}^{N}\atop {\left|I\right|=p-2}} \left(a_{I}C_{I}(z,\overline{z})+b_{I}\overline{C_{I}(z,\overline{z})}\right)+R_{2,0}(z,\overline{z}),\quad\mbox{where $R_{2,0}(z,\overline{z})\in \displaystyle\bigcap_{I\in\mathbb{N}^{N}\atop {\left|I\right|=p-2}} \ker  C^{\star}_{I}  \cap  \ker  \overline{C}^{\star}_{I} $,} \\&\quad\quad\quad\quad\vdots\quad\quad\quad\quad\quad\quad\quad\quad\vdots\quad\quad\quad\quad\quad\quad\quad\quad\vdots\quad\quad\quad\quad\quad\quad\quad\quad\vdots \quad\quad\quad\quad\quad\quad\quad\quad\quad\quad\quad\quad\quad\quad\vdots
\\& P_{k}(z,\overline{z})=P_{k+1}(z,\overline{z})q_{l}(z,\overline{z})+R_{k+1}(z,\overline{z}),\quad\mbox{where $\tr_{l}\left(R_{k+1}(z,\overline{z})\right)=0$ and:}\\& \quad\quad\quad \quad \hspace{0.1 cm} R_{k+1}(z,\overline{z})=\displaystyle\sum_{I\in\mathbb{N}^{N}\atop {\left|I\right|=p-2k}} \left(a_{I}C_{I}(z,\overline{z})+b_{I}\overline{C_{I}(z,\overline{z})}\right)+R_{k+1,0}(z,\overline{z}), \quad \mbox{where $R_{k+1,0}(z,\overline{z})\in \displaystyle\bigcap_{I\in\mathbb{N}^{N},   \atop {\left|I\right|=p-2k}} \ker  C^{\star}_{I}  \cap  \ker  \overline{C}^{\star}_{I} $,}\\&\quad\quad\quad\quad\vdots\quad\quad\quad\quad\quad\quad\quad\quad\vdots\quad\quad\quad\quad\quad\quad\quad\quad\vdots\quad\quad\quad\quad\quad\quad\quad\quad\vdots \quad\quad\quad\quad\quad\quad\quad\quad\quad\quad\quad\quad\quad\quad\vdots\end{split} \label{new1}
\end{equation}
for all $l=1,\dots,N$, where   these occurring polynomials 
\begin{equation}\left\{P_{k}(z,\overline{z})\right\}_{k=1,\dots,\left[\frac{p}{2}\right]},\quad \left\{R_{k}(z,\overline{z})\right\}_{k=1,\dots,\left[\frac{p}{2}\right]},
\label{poll1}
\end{equation}
are iteratively obtained  using the generalized version  of the Fischer Decomposition\cite{sh}.  

Recalling   strategies from \cite{bu4} and \cite{bu5}, we define 
\begin{equation}\mathcal{G}_{p}^{\left(l\right)}=\left\{\mbox{$P(z,\overline{z})$ is a real-valued polynomial of degree $p\geq 1$ in $(z,\overline{z})$ satisfying the normalizations:}\atop {P_{k}^{\left(p\right)}(z,\overline{z})=P_{k+1}^{\left(p\right)}(z,\overline{z})q_{l}(z,\overline{z})+R_{k+1}^{\left(p\right)}(z,\overline{z}),\hspace{0.1 cm}\mbox{where}\hspace{0.1 cm}R_{k+1}^{\left(p\right)}(z,\overline{z})\in\displaystyle\bigcap_{I\in\mathbb{N}^{N},\hspace{0.1 cm}I\not\in\mathcal{S}\atop {\left|I\right|=p-2k}} \ker  C^{\star}_{I}  \cap  \ker  \overline{C}^{\star}_{I}\cap \ker  \tr_{l} \hspace{0.1 cm}\mbox{and }\atop{  \mbox{for $P_{0}^{\left(p\right)}(z,\overline{z})=P(z,\overline{z})$ and   $k=0,\dots, \left[\frac{p}{2}\right]$.}} }  \right\},\quad p\in\mathbb{N}^{\star},\label{spartiuG}
\end{equation} 
for all $l=1,\dots, N$.

We consider the Fischer Decompositions (\ref{new1}) choosing 
\begin{equation}P(z,\overline{z})=\frac{\varphi_{p}(z,\overline{z})-\overline{\varphi_{p}(z,\overline{z})}}{2\sqrt{-1}},\quad\mbox{for given $p\in\mathbb{N}^{\star}$.}\label{kama1}\end{equation}

The   $G$-components of the formal transformation (\ref{map}) is  iteratively computed by  imposing the   normalizations   described by (\ref{spartiuG}). In particular, these  polynomials (\ref{poll1}) are uniquely determined by (\ref{spartiuG}), but it is not clear if these polynomials (\ref{poll1}) are real valued, because the reality of the polynomial $P(z,\overline{z})$ does not guarantee the reality of the homogeneous polynomials involved in the iterative Fischer Decompositions  described by (\ref{new1}). Then, the arguments follow \cite{val8}, because the explanations are identical. 

In order to show that these  Spaces of Normalizations (\ref{spartiuG}) uniquely determine  the $G$-component of the transformation (\ref{map}), it is required also to show the linear independence, considering complex numbers, of the following set of polynomials 
\begin{equation}\left\{C_{I,l,l'}(z,\overline{z}),\hspace{0.1 cm} \overline{C_{I,l,l'}(z,\overline{z})}\right\}_{I\in\mathbb{N}^{N}\atop {{\left|I\right|=p}  }},\quad\mbox{for all $p\in\mathbb{N}^{\star}$ and $l,l'=1,\dots,N$,}\label{330}
\end{equation}
concluding, in particular, the following
\begin{equation}a_{I,l,l'}=\overline{b_{I,l,l'}},\quad\mbox{for all $I\in\mathbb{N}^{N}$ having length $p\geq 3$ and $l,l'=1,\dots,N$.}\label{331}
\end{equation}

The non-triviality of the Fischer Decomposition\cite{sh} forces the use of (\ref{opt}) as follows: for any given multi-index 
\begin{equation}
  I=\left(i_{1},\dots,i_{k},\dots,i_{N}\right)\in\mathbb{N}^{N}\quad\mbox{such that $i_{1}+\dots+i_{k}+\dots+i_{N}=p\geq 3, $}\label{IJ1}
\end{equation}
similarly somehow to (\ref{vectori}),   we consider by (\ref{IJ1}) the following vector
\begin{equation*} Z\left[I \right]=\begin{pmatrix}a_{\left(i_{1}-2,\dots,i_{k},\dots,i_{N}\right) }^{(1)}\\ \vdots \\ a_{\left(i_{1},\dots,i_{l}-2,\dots,i_{N}\right) }^{(l)}\\ \vdots \\  a_{\left(i_{1},\dots,i_{k},\dots,i_{N}-2\right) }^{(N)}.\end{pmatrix}
\end{equation*}

Immediately from (\ref{opt}), we obtain
\begin{equation}z^{I}-A_{I,l}(z,\overline{z}) q_{l}(z,\overline{z})=C_{I,l}(z,\overline{z}),\quad\mbox{for all $I \in \mathbb{N}$ having   length $p$, for all $l=1,\dots, N$.}\label{opt1}
\end{equation} 

Clearly (\ref{opt1}) naturally introduces the following matrix
\begin{equation} 
\mbox{Aux}=  \begin{pmatrix} 1& \dots &  1& \dots & 1 \\ \vdots & \ddots & \vdots & \ddots & \vdots    \\ 1& \dots &  1 & \dots & 1\\ \vdots & \ddots & \vdots & \ddots & \vdots   \\ 1 & \dots &  1& \dots &  1
\end{pmatrix}\label{aux1},\end{equation} 
because  each very consistent  sum of terms  multiplied by $z_{l'}^{2}$ in  (\ref{opt1}), for all $l'=1,\dots,N$, generates by  (\ref{Lambda})  and (\ref{aux1})  the following terms $$\left( \mbox{Aux}\Lambda  \right)Z\left[I \right].
$$

Similarly as previously, we consider    the following matrices
 \begin{equation} \mbox{Aux}\left[I\right]= \begin{pmatrix}  1 & \dots   & 0 & \dots & 0 & \dots & 0& \dots & 0  \\ \vdots & \ddots & \vdots & \ddots & \vdots & \ddots & \vdots& \ddots & \vdots   \\ 0 & \dots & \lambda_{1} & \dots & \lambda_{k} & \dots & \lambda_{N} & \dots & 0 
\\ \vdots & \ddots & \vdots  & \ddots & \vdots  & \ddots & \vdots  & \ddots & \vdots  \\ 0 & \dots & \lambda_{1} & \dots & \lambda_{k}& \dots & \lambda_{N}& \dots & 0 \\ \vdots  & \ddots & \vdots  & \ddots & \vdots  & \ddots & \vdots  & \ddots & \vdots  \\ 0 & \dots &\lambda_{1} & \dots & \lambda_{k} & \dots & \lambda_{N}& \dots & 0 \\ \vdots  & \ddots & \vdots  & \ddots & \vdots  & \ddots & \vdots  & \ddots & \vdots  \\ 0  & \dots & 0  & \dots &0  & \dots & 0  & \dots & 1  
\end{pmatrix} .\label{calcan1}
\end{equation}

This matrix (\ref{calcan1})  has the characteristic the middle entry containing $\lambda_{k}$  is  the diagonal entry  on the  following row   
$$\left(i_{1},\dots,i_{k},\dots,i_{N}\right),$$  
according to the  lexicografic order related to (\ref{ODIN}), otherwise having  $1$  as diagonal entry,   being important in order to consider  products of matrices as in (\ref{extra4}), (\ref{extra4A}), (\ref{extra4B}), (\ref{extra4C})  in order to make  approximations of the formal transformation (\ref{map}). We define $$Z^{t}=\left\{a_{I}^{(1)},a_{I}^{(2)},\dots,a_{I}^{\left(N\right)},\right\}_{I\in\mathbb{N}^{N}\atop \left|I\right|=p},$$
which respects the lexicografic order related to (\ref{ODIN}).  

We obtain the following system    of  equations 
\begin{equation}\left(I_{N^{p}}-\mbox{Aux}_{p}A\right) Z+B\overline{Z}=V\left(z_{1},z_{2},\dots,z_{N}\right),\quad\mbox{where $A$, $B\in \mathcal{M}_{N^{p}\times N^{p}}\left(\mathbb{C}\right)$}, \label{sisi}
\end{equation} 
and $V\left(z_{1},z_{2},\dots,z_{N}\right)$ is a known vector homogeneous polynomial  of degree $p-2$, dealing with the following products of ,,simple'' matrices
\begin{equation} 
\mbox{Aux}_{p}=\mbox{Aux}\left[p-2,\dots,0\right]\mbox{Aux}\left[p-3,1,\dots,0\right]\dots\mbox{Aux}\left[0,\dots,1,p-3\right] \mbox{Aux}\left[0,\dots,p-2\right],\quad\mbox{for all $p \geq 3$.}\label{calcan2}
\end{equation}

\subsection{Fischer Normalization $F$-Spaces\cite{bu2},\cite{bu3}}The Fischer Decomposition from (\ref{opsec}) gives
\begin{equation} \left(\overline{z}_{l}+2\lambda_{l}z_{l}\right) z^{J}=A_{l,J,l'}(z,\overline{z}) q_{l'}(z,\overline{z})+C_{l,J,l'}(z,\overline{z}),\quad  \tr \left(C_{l,J,l'}(z,\overline{z})\right)=0,\quad\mbox{where      $l',l= 1,\dots, N$.} \label{optsec} 
\end{equation}
 
This set of polynomials, derived from (\ref{optsec}),   dictates   normalization conditions defining  certain  Spaces of Fischer Normalizations, which are constructed   iteratively from the generalized version of  the Fischer Decomposition\cite{sh}: for a given  real homogeneous polynomial of degree $p\geq 1$ in $(z,\overline{z})$ denoted by $P(z,\overline{z})$,  we  have
\begin{equation} \begin{split}& P(z,\overline{z})=P_{1}(z,\overline{z})q_{l'}(z,\overline{z})+R_{1}(z,\overline{z}),\quad\mbox{where $\tr_{l'}\left(R_{1}(z,\overline{z})\right)=0$ such that:}\\& \quad  R_{1}(z,\overline{z})=\displaystyle\sum_{l=1}^{N}\displaystyle\sum_{J\in\mathbb{N}^{N}\atop {\left|J\right|=p-1}} \left(a_{l,J}C_{l,J}(z,\overline{z})+b_{l,J}\overline{C_{l,J}(z,\overline{z})}\right)+R_{1,0}(z,\overline{z}),\quad\mbox{for:}  \quad \mbox{ $R_{1,0}(z,\overline{z})\in \displaystyle\bigcap_{l=1}^{N}\displaystyle\bigcap_{J\in\mathbb{N}^{N} \atop {\left|J\right|=p-1}} \ker  C^{\star}_{l,J}  \bigcap  \ker  \overline{C}^{\star}_{l,J} $,}\\& P_{1}(z,\overline{z})=P_{2}(z,\overline{z})q_{l'}(z,\overline{z})+R_{2}(z,\overline{z}),\quad\mbox{for: $ \tr_{l'}\left(R_{2}(z,\overline{z})\right)=0$ such that:}\\& \quad    R_{2}(z,\overline{z})=\displaystyle\sum_{l=1}^{N}\displaystyle\sum_{J\in\mathbb{N}^{N}\atop {\left|J\right|=p-3}} \left(a_{l,J}C_{l,J}(z,\overline{z})+b_{l,J}\overline{C_{l,J}(z,\overline{z})}\right)+R_{2,0}(z,\overline{z}),\quad\mbox{where:}   \quad  \mbox{ $R_{2,0}(z,\overline{z})\in \displaystyle\bigcap_{l=1}^{N}\displaystyle\bigcap_{J\in\mathbb{N}^{N}, \atop {\left|J\right|=p-3}} \ker  C^{\star}_{l,J}  \bigcap   \ker  \overline{C}^{\star}_{l,J} $,}
\\&\quad\quad\quad\quad\vdots\quad\quad\quad\quad\quad\quad\quad\quad\vdots\quad\quad\quad\quad\quad\quad\quad\quad\vdots\quad\quad\quad\quad\quad\quad\quad\quad\vdots \quad\quad\quad\quad\quad\quad\quad\quad\quad\quad\quad\quad\quad\quad\vdots\\& P_{k}(z,\overline{z})=P_{k+1}(z,\overline{z})q_{l}(z,\overline{z})+R_{k+1}(z,\overline{z}),\quad\mbox{where $\tr_{l'}\left(R_{k+1}(z,\overline{z})\right)=0$ such that:}\\&   R_{k+1}(z,\overline{z})=\displaystyle\sum_{l=1}^{N}\displaystyle\sum_{J\in\mathbb{N}^{N}\atop {\left|J\right|=p-2k-1}} \left(a_{l,J}C_{l,J}(z,\overline{z})+b_{l,J}\overline{C_{l,J}(z,\overline{z})}\right)+R_{k+1,0}(z,\overline{z}),\quad\mbox{for:}   \hspace{0.1 cm}  \mbox{ $R_{k+1,0}(z,\overline{z})\in \displaystyle\bigcap_{l=1}^{N}\displaystyle\bigcap_{J\in\mathbb{N}^{N}, \atop {\left|J\right|=p-2k-1}} \ker  C^{\star}_{l,J}  \bigcap   \ker  \overline{C}^{\star}_{l,J} $,}\\&\quad\quad\quad\vdots\quad\quad\quad\quad\quad\quad\quad\quad\vdots\quad\quad\quad\quad\quad\quad\quad\quad\vdots\quad\quad\quad\quad\quad\quad\quad\quad\vdots \quad\quad\quad\quad\quad\quad\quad\quad\quad\quad\quad\quad\quad\quad\vdots\end{split}.\label{new2}
\end{equation}
for all $l'=1,\dots,N$, where   these occurring polynomials 
\begin{equation}\left\{P_{k}(z,\overline{z})\right\}_{k=1,\dots,\left[\frac{p-1}{2}\right]},\quad \left\{R_{k}(z,\overline{z})\right\}_{k=1,\dots,\left[\frac{p-1}{2}\right]},
\label{poll2}
\end{equation}
are iteratively obtained  using the generalized version  of the Fischer Decomposition\cite{sh}. 

Recalling   strategies from \cite{bu2} and \cite{bu3}, we define 
\begin{equation}\mathcal{F}_{p}^{\left(l'\right)}=\left\{\mbox{$P(z,\overline{z})$ is a real-valued polynomial of degree $p\geq 1$ in $(z,\overline{z})$ satisfying the normalizations:}\atop {P_{k}^{\left(p\right)}(z,\overline{z})=P_{k+1}^{\left(p\right)}(z,\overline{z})q_{l'}(z,\overline{z})+R_{k+1}^{\left(p\right)}(z,\overline{z}),\hspace{0.1 cm}\mbox{where}\hspace{0.1 cm}R_{k+1}^{\left(p\right)}(z,\overline{z})\in \displaystyle\bigcap_{l=1}^{N}\displaystyle\bigcap_{J\in\mathbb{N}^{N}\atop {\left|J\right|=p-2k-1}} \ker  C^{\star}_{l,J}  \bigcap   \ker  \overline{C}^{\star}_{l,J}\bigcap  \ker\tr_{l'}, \hspace{0.1 cm}\mbox{and}\atop{\mbox{ \mbox{for}\hspace{0.1 cm} $P_{0}^{\left(p\right)}(z,\overline{z})=P(z,\overline{z})$ and    $k=0,\dots, \left[\frac{p-1}{2}\right]$ .}} }  \right\},\quad p\in\mathbb{N}^{\star},\label{spartiuF}
\end{equation}
for all $l'=1,\dots,N$ 
 
We  consider the Fischer Decompositions (\ref{new1}) choosing 
\begin{equation}P(z,\overline{z})=\frac{\varphi_{p}(z,\overline{z})+\overline{\varphi_{p}(z,\overline{z})}}{2},\quad\mbox{for given $p\in\mathbb{N}^{\star}$.}\label{kama2}\end{equation}

In order to show now that these Spaces of Normalizations (\ref{spartiuF}) uniquely determine  the $F$-component of the formal transformation (\ref{map}), it is required to show the linear independence, considering complex numbers, of the following set of polynomials 
\begin{equation}\left\{C_{l,J,l'}(z,\overline{z}),\hspace{0.1 cm} \overline{C_{l,J,l'}(z,\overline{z})}\right\}_{J\in\mathbb{N}^{N}\atop{\left|J\right|=p-1 \atop{l,l'=1,\dots,N}}},\quad\mbox{for all $p\in\mathbb{N}^{\star}$.}\label{Africa330}
\end{equation}

In particular, showing the linear independence of these polynomials   implies
\begin{equation}a_{l,J,l'}=\overline{b_{l,J,l'}},\quad\mbox{for all $J\in\mathbb{N}^{N}$ having length $p-1\geq 2$ and  $l,l'=1,\dots,N$.}\label{Africa331}
\end{equation}

These computations are difficult  to conclude because of  the overlapping of the  homogeneous   polynomials from (\ref{Africa330}). There are recalled   the previous arguments  from \cite{val8}.  It is necessary   to study more carefully all their corresponding computations of its solution as follows: for any given multi-index 
\begin{equation}
  J=\left(j_{1},\dots,j_{k},\dots,j_{N}\right)\in\mathbb{N}^{N}\quad\mbox{such that $j_{1}+\dots+j_{k}+\dots+j_{N}=p-1\geq 2$,}\label{AfricaIJ1}
\end{equation}
similarly somehow to (\ref{vectori}), we consider by (\ref{AfricaIJ1}) the following vector
\begin{equation} \tilde{Z}\left[J \right]=\begin{pmatrix}a_{\left(j_{1}-2,\dots,j_{k},\dots,j_{N}\right)}^{\left(1\right)}\\ \vdots \\ a_{\left(j_{1},\dots,j_{k}-2,\dots,j_{N}\right)}^{\left(1\right)}\\ \vdots \\  a_{\left(j_{1},\dots,j_{k},\dots,j_{N}-2\right)}^{\left(1\right)}\\  \vdots\\ a_{\left(j_{1}-2,\dots,j_{k},\dots,j_{N}\right)}^{(N)}\\ \vdots \\ a_{\left(j_{1},\dots,j_{k}-2,\dots,j_{N}\right)}^{(N)}\\ \vdots \\  a_{\left(j_{1},\dots,j_{k},\dots,j_{N}-2\right)}^{(N)}.\end{pmatrix}\label{kaka}
\end{equation}

Immediately from (\ref{opt}), we obtain
\begin{equation}\left(\overline{z}_{l}+2\lambda_{l}z_{l}\right) z^{J}-A_{l,J,l'}(z,\overline{z}) q_{l'}(z,\overline{z})=C_{l,J,l'}(z,\overline{z}),\quad  \tr_{l'}\left(C_{l,J,l'}(z,\overline{z})\right)=0,\quad\mbox{for all  $l,l'=1,\dots, N$.}\label{Africaopt1}
\end{equation} 

Then each very consistent  sum of terms  multiplied by $\lambda_{l}z_{l}^{2}$, for all $l=1,\dots, N$ in  (\ref{optsec}), generates by  (\ref{Lambda}), (\ref{aux1}) and (\ref{AfricaIJ1}) the following terms\begin{equation} 
 \begin{pmatrix} \mbox{Aux}\Lambda  &  \mbox{O}_{N^{p-1}}& \dots & \mbox{O}_{N^{p-1}}     \\ \mbox{O}_{N^{p-1}}&    \mbox{Aux}\Lambda  & \dots & \mbox{O}_{N^{p-1}}\\ \vdots   & \vdots & \ddots & \vdots   \\ \mbox{O}_{N^{p-1}}   &  \mbox{O}_{N^{p-1}}& \dots &  \mbox{Aux}\Lambda 
\end{pmatrix}\tilde{Z}\left[J \right].\label{Africaaux1}\end{equation} 

We consider  by (\ref{AfricaIJ1}) the following matrices
 \begin{equation} \tilde{\mbox{Aux}}\left[J\right]=  \begin{pmatrix} \mbox{Aux}\left[J\right]  &  \mbox{O}_{N^{p-1}}& \dots & \mbox{O}_{N^{p-1}}     \\ \mbox{O}_{N^{p-1}}&    \mbox{Aux}\left[J\right] & \dots & \mbox{O}_{N^{p-1}}\\ \vdots   & \vdots & \ddots & \vdots   \\ \mbox{O}_{N^{p-1}}   &  \mbox{O}_{N^{p-1}}& \dots &  \mbox{Aux}\left[J\right] 
\end{pmatrix} .\label{Africacalcan1}
\end{equation}

 Defining $$\tilde{Z}^{t}=\left(\left\{a_{J}^{\left(1\right)}\right\}_{J\in\mathbb{N}^{N}\atop \left|I\right|=p-1},\left\{a_{J}^{(2)}\right\}_{J\in\mathbb{N}^{N}\atop \left|I\right|=p-1},\dots,\left\{a_{J}^{(N)}\right\}_{J\in\mathbb{N}^{N}\atop \left|I\right|=p-1}\right),$$
according to the lexicografic order related to (\ref{ODIN}), we obtain the following  system   of  equations 
\begin{equation}\left(I-\mbox{Aux}_{p}\tilde{A}\right) \tilde{Z}+\tilde{B}\overline{\tilde{Z}}=\tilde{V}\left(z_{1},z_{2},\dots,z_{N}\right), \label{Africasisise}
\end{equation} 
where $\tilde{V}\left(z_{1},z_{2},\dots,z_{N}\right)$ is a known homogeneous vector polynomial  of degree $p-2$, dealing with the following products of ,,simple'' matrices 
\begin{equation} 
\tilde{\mbox{Aux}}_{p}=\tilde{\mbox{Aux}}\left[p-2,\dots,0\right]\tilde{\mbox{Aux}}\left[p-3,1,\dots,0\right]\dots\tilde{\mbox{Aux}}\left[0,\dots,1,p-3\right] \tilde{\mbox{Aux}}\left[0,\dots,p-2\right],\quad\mbox{for all $p \geq 3$.}\label{Africacalcan2}
\end{equation}

We were using the following notations
\begin{equation}\tilde{A}=  \begin{pmatrix} A_{1} &  \mbox{O}_{N^{p-1}}& \dots & \mbox{O}_{N^{p-1}}     \\ \mbox{O}_{N^{p-1}}&    A_{2} & \dots & \mbox{O}_{N^{p-1}}\\ \vdots   & \vdots & \ddots & \vdots   \\ \mbox{O}_{N^{p-1}}   &  \mbox{O}_{N^{p-1}}& \dots &  A_{N}
\end{pmatrix} ,\quad \tilde{B}=  \begin{pmatrix} B_{1} &  \mbox{O}_{N^{p-1}}& \dots & \mbox{O}_{N^{p-1}}     \\ \mbox{O}_{N^{p-1}}&    B_{2} & \dots & \mbox{O}_{N^{p-1}}\\ \vdots   & \vdots & \ddots & \vdots   \\ \mbox{O}_{N^{p-1}}   &  \mbox{O}_{N^{p-1}}& \dots &  B_{N}
\end{pmatrix} ,\label{Africacalcan1}
\end{equation}
where   we have
$$A_{1}, A_{2},\dots, A_{N}; B_{1},  B_{2}, \dots,  B_{N}\in \mathcal{M}_{N^{p-1}\times N^{p-1} }\left(\mathbb{C}\right).
$$

It is   clear that (\ref{Africasisise}) has unique solution due to the existences of the following matrices
\begin{equation}\frac{1}{I_{N^{p}}-\tilde{\mbox{Aux}}_{p}\tilde{A}-\tilde{B}},\quad\frac{1}{I_{N^{p}}-\tilde{\mbox{Aux}}_{p}\tilde{A}+\tilde{B}},\label{5501}
\end{equation}
which follow from \cite{val8}.

\section{Proofs of Theorems  \ref{ta}  and  \ref{tB} }
\subsection{Proof of Theorem  \ref{ta} } Clearly, any such formal embedding may be written as (\ref{map}). Then, up to compositions with suitable automorphisms of corresponding models, we obtain just a class of equivalence, which results from the uniqueness of the Fischer Decompositions related to (\ref{spartiuG}) and (\ref{spartiuF}). 
\subsection{Proof of Theorem  \ref{tB} } It follows easily as in \cite{val8}.

\end{document}